\documentclass[article,12pt]{article} 
\usepackage{amscd,amsmath,amssymb,amsfonts,latexsym,mathrsfs,amsthm}
\usepackage{graphicx} 
\usepackage{setspace} 
\usepackage{graphics,epsfig}
\usepackage{sectsty}
\usepackage{natbib}
\usepackage[english]{babel}
\usepackage{amsmath}
\usepackage{amsfonts}
\usepackage{amssymb,amsthm}
\usepackage{bm}
\usepackage{mathrsfs}
\usepackage{subfigure}
\usepackage[usenames]{color}
\usepackage{rotating}

\usepackage{url}

\newtheorem{rem}{{\bf Remark}}

\usepackage{flushend}
\usepackage{verbatim}
\usepackage{tabularx}
\usepackage{multirow} 
\usepackage{arydshln}

\UseRawInputEncoding

\title{Effective sample size approximations as entropy measures} 

\author{L. Martino$^\star$, V. Elvira$^\top$, \\
{\small$^\star$  Universit{\'a} degli Studi di Catania, Italy.}\\
{\small $^\top$ University of Edinburgh, UK.} \\
}
\date{}
 
\begin{document}

\maketitle

\thispagestyle{empty}

\begin{abstract}
 In this work, we analyze alternative effective sample size (ESS) metrics for importance sampling algorithms, and discuss a possible extended range of applications.  We show the relationship between the ESS expressions used in the literature and two entropy families, the R\'enyi and Tsallis entropy. The R\'enyi entropy is connected to the Huggins-Roy's ESS family introduced in \cite{Huggins15}. We prove that that all the ESS functions included in the Huggins-Roy's family fulfill all the desirable theoretical conditions.  
We analyzed and remark the connections with several other fields, such as the Hill numbers introduced in ecology, the Gini inequality coefficient employed in economics,  and the Gini impurity index used mainly in machine learning, to name a few.   
  Finally, by numerical simulations, we study the performance of different ESS expressions contained in the previous ESS families in terms of approximation of the theoretical ESS definition, and show the application of ESS formulas in a variable selection problem. 
\newline
\newline 
{ \bf Keywords:} 
Effective Sample Size; Importance Sampling; Entropy; Diversity measure; Gini impurity; Gini inequality coefficient; inverse Simpson concentration; Berger-Parker index.
\end{abstract} 

\begin{figure}[!hbt]
\centering 
\centerline{
\includegraphics[width=12cm]{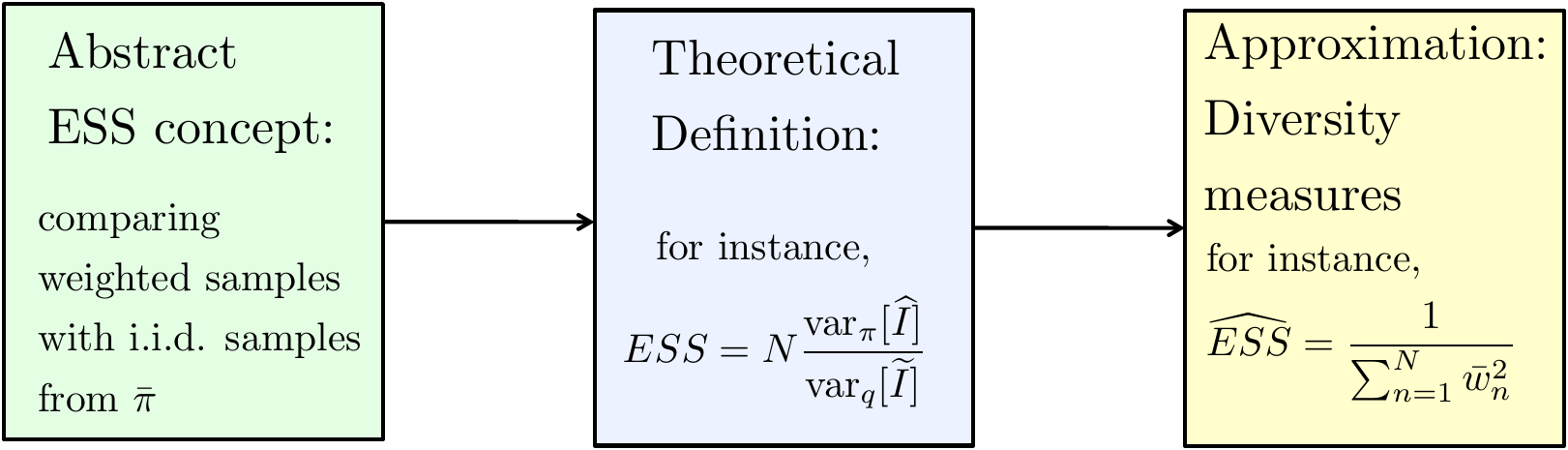}
}
  \caption{{\small Graphical representation of the development of the approximated ESS formulas for importance sampling. The abstract concept of Effective Sample Size has been translated in a mathematical formulation providing a first attempt of theoretical definition. Since this definition cannot compute, several approximations have been proposed  (based only in the information provided by the normalized IS weights). The expression $\mbox{ESS}=\frac{1}{\sum_{n=1}^M {\bar w}_n^2}$ is the most applied so far in the literature.  } }
\label{figESSteo}
\end{figure}

The effective sample size (ESS) measure is an important concept in order to quantify the efficiency of different Monte Carlo methods, such as Markov Chain Monte Carlo (MCMC) \citep{GameRman_Book,Liang10} and Importance Sampling (IS) techniques \citep{Bugallo15,Cappe04}. 
In an IS context, the ESS is a heuristic to approximate
how many independent identically distributed (i.i.d.) samples, drawn directly from the target distribution $\bar{\pi}({\bf x})=\frac{1}{Z} \pi({\bf x})$ where $Z$ is the normalizing constant, are equivalent {\it in some sense} to the $N$ weighted samples, ${\bf x}_1,\ldots, {\bf x}_n$, drawn from a proposal distribution $q({\bf x})$  and weighted according
to the ratio $w_n=\frac{\pi({\bf x}_n)}{q({\bf x}_n)}$ \citep{Robert10b}.  This consideration is represented in the first box of Figure \ref{figESSteo}, referred as {\em abstract ESS concept}.
\newline
The theoretical definition of the ESS for IS is given by the ratio between two variances \citep{GameRman_Book,Kong92}: the variance of the ideal Monte Carlo estimator (drawing samples directly from the target), and the variance of the estimator obtained by an IS scheme, using the same number of samples in both estimators (see Eq. \eqref{DEF_ESS_1} for more details). This definition presents some drawbacks (see \citep{ESSarxiv16,elvira2018rethinking} for an exhaustive discussion) and is useless for practical purposes since it cannot be computed in general. Hence, approximations of this theoretical formula are required.
In Figure \ref{figESSteo}, this theoretical definition is represented by the second box. Within an IS context, the most common choice in literature to approximate this theoretical ESS definition is  $\mbox{ESS}=\frac{1}{\sum_{n=1}^N {\bar w}_n^2}$, which involves (only) the normalized importance weights ${\bar w}_n=\frac{w_n}{\sum_{j=1}^N w_j}$, $n=1,\ldots,N$ \citep{Djuric03,SMC01,Kong94,Robert10b}. 
This expression has been  widely used in particle filtering in order to apply the resampling steps adaptively \citep{SMC01,Djuric03,Gordon93}.  However, it presents different weaknesses since it has been obtained after several approximations of the theoretical definition. For instance, it just depends on the normalized weights, but it is not dependent on particle locations and from the particular integral to approximate (see \citep{elvira2018rethinking,ESSarxiv16} for further details). Several other alternatives have been studied in literature and applied in order to perform adaptive resampling within sequential Monte Carlo (SMC) methods \cite{Huggins15,ESSarxiv16}. For instance, another measure called  perplexity, involving the discrete entropy \citep{Cover91} of the normalized weights has been also proposed in \citep{pmc-cappe08} (see also \cite[Chapter 4]{Robert10b}, \cite[Section 3.5]{Doucet08tut}). Another expression is defined as the inverse of the maximum of the normalized weights ${\bar w}_n$ \citep{ESSarxiv16}.
\newline
\newline
In this work, %we explain why these formulas have been applied so widely even 
we recall the definition of the generalized ESS (G-ESS) functions  given in \citep{ESSarxiv16}. We stress and show that the G-ESS functions can be considered {\it diversity indices} \citep{Jost06} (see third box in Figure \ref{figESSteo}). %The actual reason for the success of the ESS expressions introduced in the literature, is related to the fact they are discrepancy measure and/or can considered as diversity indices.
 Indeed, we show that the G-ESS functions can be associated to different entropy families \citep{Cover91}. Given an entropy measure of the probability mass function (pmf) defined by the normalized weights ${\bar w}_n$, $n=1,\ldots,N$,  we can obtain a G-ESS formula by taking the exponential transformation of the entropy expression (in some cases, some additional translation and scaling are needed). 
More specifically, we analyze the R\'enyi and Tsallis entropy families, converting them in G-ESS functions. The  ESS formulas corresponding to the R\'enyi entropy coincides with the Huggins-Roy's ESS family introduced and studied independently in \citep{Huggins15},
\begin{eqnarray*}
\mbox{ESS}
&=&\left(\sum_{n=1}^N {\bar w}_n^\beta \right)^{\frac{1}{1-\beta}}, Ê\quad \quad \beta\geq 0.
\end{eqnarray*}
 We show that all the G-ESS expressions belonging to this family satisfy all the desired requirements, being all defined as {\it proper and stable} (see Section \ref{GESSsect} for further details).
Moreover, almost all the main formulas previously proposed in the literature are contained in the Huggins-Roy's family. Using the Tsallis entropy, we obtain another ESS family which contain the {\it Gini impurity index} as special case, that is widely employed in  machine learning within decision tree algorithms \cite{bishop2007,Krzywinski2017}. We also discuss the connection to another ESS family provided in \citep{ESSarxiv16}.  However, generally the Tsallis ESS formulas are not proper and stable.  
 Other stable expressions that do not belong  to the Huggins-Roy's family are also given (see, e.g., Sections \ref{other_ESS_Stable} and \ref{ESSpolSect}).
\newline
 The connections with the entropy families show the relationships with multiple studies in different fields (e.g., ecology and machine learning to name a few).  The benefit  of creating these {\it bridges} between fields is twofold and bidirectional: different ideas used in other fields can be applied as ESS expressions in an IS context (as the formula \eqref{ESSGOL} introduced in political science) and, vice-versa, ESS formulas proposed for IS could be employed in other fields. Showing these bridges is the main goal of this work. 
We remark on the links with other fields in Section \ref{EcoSect0}, where we discuss the applications of ESS expression in ecology, economics, political science, physics, and also in feature selection problems.  Other connections with machine learning, economics, and ecology are also discussed in the previous Sections \ref{EntropiesSect} and \ref{other_ESS_Stable}. Figure \ref{allNames} provides a summary of the main nomenclature employed in different fields.
\newline
Furthermore, by numerical simulations, we obtain the G-ESS function within Huggins-Roy's family which provides the best approximation the theoretical ESS definition, in two specific scenarios. We also study linear combinations of G-ESS functions in order to enhance the approximation of the theoretical definition. The results of our numerical simulations suggest  the use of the formulas of the type of $\mbox{ESS}=\left(\frac{1}{\sum_{n=1}^N {\bar w}_n^4}\right)^{1/3}$ and  $\mbox{ESS}=\left(\frac{1}{\sum_{n=1}^N {\bar w}_n^{8}}\right)^{1/7}$.  Both expression differ from the classical formula $\mbox{ESS}=\frac{1}{\sum_{n=1}^N {\bar w}_n^2}$, which is contained in Huggins-Roy's family with $\beta=2$. Our study suggest the use of $\beta >2$. Moreover, we have applied the most relevant ESS formulas in a variable selection framework. Some of them provides good results in line with the expert's opinions.
  These considerations can be also relevant clues for future applications and studies.

 %%%%%%%%%%%%%%%%%%%%%%%

%%%%%%%%%%%%%%%%%%%%%%
 %{\bf Keywords:}
 %%%%%%%%%%%%%%%%%%%%%%%
 
 %%%%%%%%%%%%%%%%%%%%%%% %%%%%%%%%%%%%%%%%%%%%%% %%%%%%%%%%%%%%%%%%%%%%%
\section{Effective sample size (ESS) for importance sampling}
 %%%%%%%%%%%%%%%%%%%%%%% %%%%%%%%%%%%%%%%%%%%%%% %%%%%%%%%%%%%%%%%%%%%%%%
 \label{ClaDer}
 Let us denote the target probability density function (pdf)  as $\bar{\pi}({\bf x}) \propto \pi({\bf x})$ (known up to a normalizing constant) with ${\bf x}\in\mathcal{X}$. Moreover, we consider the following integral involving $\bar{\pi}({\bf x})$ and a square-integrable function $h({\bf x})$,   
 \begin{equation}
I=\int_{\mathcal{X}} h({\bf x}) \bar{\pi}({\bf x}) d{\bf x},
\end{equation}
which we aim to approximate using a  Monte Carlo approach. If we are able to draw $N$ independent samples ${\bf x}_1,\ldots,{\bf x}_N$ from $\bar{\pi}({\bf x})$, then the Monte Carlo estimator of $I$ is 
 \begin{equation}\label{IdealMCestimator}
{\widehat I}=\frac{1}{N}\sum_{n=1}^N   h({\bf x}_n)%\xrightarrow[\text{world}]{\text{hello}}
\xrightarrow[]{N \rightarrow \infty}
 I,  \quad \mbox{ where } \quad {\bf x}_n \sim  \bar{\pi}({\bf x}).
\end{equation}
 However, generating samples directly from the target,  $\bar{\pi}({\bf x})$, is often impossible. Alternatively, we can draw  $N$  samples ${\bf x}_1,\ldots,{\bf x}_N$ from a (simpler) proposal pdf $q({\bf x})$,\footnote{We assume that $q({\bf x})>0$ for all ${\bf x}$ where $\bar{\pi}({\bf x})) \neq 0$.} and then assign a weight  to each sample,
$ w_n =\frac{\pi({\bf x}_n)}{q({\bf x}_n)}$, with $n=1,\ldots,N$, according to the importance sampling (IS) approach. Defining the normalized weights,
\begin{equation}
\label{NormWeights}
{\bar w}_n=\frac{ w_n }{\sum_{i=1}^N w_i}, \quad \mbox{ where } \quad  w_n =\frac{\pi({\bf x}_n)}{q({\bf x}_n)},  \quad n=1,\ldots,N,
\end{equation}
then the self-normalized IS estimator is 
 \begin{equation}
{\widetilde I}=\sum_{n=1}^N  {\bar w}_n h({\bf x}_n)\xrightarrow[]{N \rightarrow \infty} I, \quad \mbox{ where } \quad {\bf x}_n \sim q({\bf x}).
\end{equation}
Generally, the estimator ${\widetilde I}$ has greater variance than ${\widehat I}$, since the samples are not directly drawn from $\bar{\pi}({\bf x})$ (for some exceptions, that occur with a suitable choice of the proposal, see \cite{optimalityIS}).  Moreover,  ${\widetilde I}$ is biased whereas ${\widehat I}$ is unbiased.
 In several applications \citep{Djuric03,SMC01}, it is necessary to measure the loss of the efficiency when we apply the IS estimator ${\widetilde I}$, instead of ideal Monte Carlo estimator ${\widehat I}$, i.e., to measure in some way the increase of variance due to the use of ${\widetilde I}$ instead of  ${\widehat I}$. Hence, the idea is to define the Effective Sample Size (ESS) as the ratio of the variances of the estimators \citep{Kong92},
\begin{equation}
\label{DEF_ESS_1}
\mbox{ESS}_{\texttt{teo}}(h)=N\frac{\mbox{var}_\pi[{\widehat I}]}{\mbox{var}_q[{\widetilde I}]}.
\end{equation}
Note the dependence on the function $h({\bf x})$ corresponding to a specific integral.

 %%%%%%%%%%%%%%%%%%%%%%%%%%%%%%%%%%%
\section{Practical ESS formulas}
 %%%%%%%%%%%%%%%%%%%%%%% %%%%%%
 %%%%%%%%%%%%%%%%%%%%%%%%%%%%%%%%%%%
\subsection{ESS expressions in the literature}
 %%%%%%%%%%%%%%%%%%%%%%% %%%%%%
Finding a useful expression of ESS derived analytically from the theoretical definition in Eq. \eqref{DEF_ESS_1} above is not straightforward. Then, different derivations \citep{Kong92,Kong94}, \cite[Chapter 11]{SMC01}, \cite[Chapter 4]{Robert10b} proceed using several approximations and assumptions for yielding an expression useful from a practical point of view.  A  well-known rule of thumb, widely used in literature \citep{SMC01,Liu04b,Robert10b}, is 
\begin{eqnarray}
\label{FirstDef_P}
\mbox{ESS}_N({\bf {\bar w}})&=& \frac{1}{\sum_{n=1}^N {\bar w}_n^2}, 
\end{eqnarray} 
where we have used the the normalized weights 
$$
{\bf {\bar w}}=[{\bar w}_1, \ldots, {\bar w}_N],
$$ 
defined in Eq. \eqref{NormWeights}. The formula above 
 has also an intuitive probabilistic interpretation (from a resampling point of view): if we draw random pairs of samples with replacement according to the probability mass function (pmf) defined  by ${\bar w}_n$, with $n=1,...,N$, the value  $ \frac{1}{\sum_{n=1}^N {\bar w}_n^2}$ is the expected number of trials needed to obtain a {\it first} pair containing the same sample twice (see Appendix \ref{App1} for details).  Furthermore, another interesting form of Eq. \eqref{FirstDef_P} as a function of the variance of the weights is given in Appendix \ref{App2}.  Another similar measure, called {\it perplexity},  has been proposed independently in literature  based only  on the normalized importance weights \citep{pmc-cappe08,Robert10b},
\begin{eqnarray}
\label{PerpEq}
\mbox{ESS}_N({\bf {\bar w}})=\exp\{{H({\bf {\bar w}})}\}
\end{eqnarray}  
where 
$$
H({\bf {\bar w}})= - \sum_{n=1}^N \bar{w}_n \log \bar{w}_n,
$$
is the discrete entropy of the vector ${\bf {\bar w}}$ \citep{Cover91}.
An additional example is the following formula \citep{ESSarxiv16},
\begin{eqnarray}
\label{MaxEq}  
\mbox{ESS}_N({\bf {\bar w}})=\frac{1}{\max \bar{w}_n}.
\end{eqnarray} 
 Let us assume that $\max \bar{w}_n$  is reached only with one sample (only for one index $n$). In this case, the expression above 
 has also a probabilistic interpretation: if we draw one sample with replacement according to the pmf defined  by ${\bar w}_n$, with $n=1,...,N$, the value  $\frac{1}{\max \bar{w}_n}$ is the expected number of trials needed to obtain for a first time  the sample corresponding to the maximum weight. The proof is very similar the derivation in Appendix \ref{App1}.  This interpretation is interesting from a resampling point of view. An interesting property of all the three expressions above in Eqs. \eqref{FirstDef_P}-\eqref{PerpEq}-\eqref{MaxEq} is
\begin{equation}\label{Constrains_ESS}
1 \leq \mbox{ESS}_N({\bf {\bar w}})\leq N.
\end{equation}

%%%%%%%%%%%%%%%%%%%%%%%%%%%%%%%%
 \subsection{Relationship with the theoretical definition}
 %%%%%%%%%%%%%%%%%%%%%%%%%%%%%%%%

All these measures $\mbox{ESS}_N({\bf {\bar w}})$ are only based on the normalized weights ${\bf {\bar w}}$ and  there is a loss of   information regarding the locations of the samples ${\bf x}_n$, which is clearly a drawback \cite{elvira2018rethinking,ESSarxiv16}, 
even if the computation of the weights involves the use of the samples, i.e.,  $ w_n =\pi({\bf x}_n)/q({\bf x}_n)$.
 To  clarify this point, we give the following example.  Two different samples $\mathbf{x}'$ and $\mathbf{x}''$ could have  very similar weights $w'=\frac{\pi\left(\mathbf{x}'\right)}{q\left(\mathbf{x}'\right)}\approx w''=\frac{\pi\left(\mathbf{x}''\right)}{q\left(\mathbf{x}''\right)}$, so that the ESS formulas just use this information. However, the ESS formulas lose all the information about the positions of the samples $\mathbf{x}'$ and $\mathbf{x}''$. The two particles can be very close to each other or far away; the latter scenario is often preferred in terms of statistical information. 
\newline
Moreover, the theoretical value $\mbox{ESS}_{\texttt{teo}}(h)$ in \eqref{DEF_ESS_1} is always positive, could be smaller than $1$ and, in some situations, bigger than $N$ as well  \cite[Section 3.3]{elvira2018rethinking}, \cite{optimalityIS}. This last scenario can occur 
 when an optimal proposal pdf (or a density close to the optimal one) is used in an IS scheme \cite{optimalityIS}. In this case, the IS scheme can beat the baseline Monte Carlo estimator, and 
$\mbox{ESS}_{\texttt{teo}}(h) >N$.  Furthermore, $\mbox{ESS}_{\texttt{teo}}(h)$ depends on the function $h$ that does not appear in the expressions $\mbox{ESS}_N({\bf {\bar w}})$.
Therefore, the formulas $\mbox{ESS}_N({\bf {\bar w}})$ that all satisfy the constrains in Eq. \eqref{Constrains_ESS} (i.e., $1 \leq \mbox{ESS}_N({\bf {\bar w}})\leq N$) are quite rough approximations of   $\mbox{ESS}_{\texttt{teo}}(h)$. However, they are often used in practice. The reason for this success is connected to their interpretation as discrepancy/diversity measures, as explained below.

%%%%%%%%%%%%%%%%%%%%%%%%%%%%
 \subsection{Discrepancy w.r.t. the uniform pmf}
  %%%%%%%%%%%%%%%%%%%%%%%%%%
  All the formulas above can be considered {\it diversity indices} or {\it discrepancy measures} \cite{Jost06,ESSarxiv16}. We give more details in the rest of the work. Here,  let us start considering the discrepancy between two pmfs: the pmf defined by the weights $\mathbf{{\bar w}}=\left[\bar{w}_1, \ldots, \bar{w}_N\right]$ and the discrete uniform pmf defined by $\mathbf{{\bar w}}^*=\left[\frac{1}{N}, \ldots, \frac{1}{N}\right]$. Indeed, the ESS formula in Eq. \eqref{FirstDef_P} can be directly related to the Euclidean distance between these two pmfs, i.e.,
$$
\begin{aligned}
\left\|\mathbf{\bar{w}}-\mathbf{{\bar w}}^*\right\|_2 & =\sqrt{\sum_{n=1}^N\left(\bar{w}_n-\frac{1}{N}\right)^2} \\
& =\sqrt{\left(\sum_{n=1}^N \bar{w}_n^2\right)+N\left(\frac{1}{N^2}\right)-\frac{2}{N} \sum_{n=1}^N \bar{w}_n} \\
& =\sqrt{\left(\sum_{n=1}^N \bar{w}_n^2\right)-\frac{1}{N}} \\
& =\sqrt{\frac{1}{\mbox{ESS}_N({\bf {\bar w}})}-\frac{1}{N}} ,
\end{aligned}
$$
where we have used $\mbox{ESS}_N({\bf {\bar w}})=\frac{1}{\sum_{i=1}^N {\bar w}_n^2}$  in Eq. \eqref{FirstDef_P}.
Hence, maximizing the expression in Eq. \eqref{FirstDef_P} is equivalent to minimizing the Euclidean distance $\left\|\mathbf{\bar{w}}-\mathbf{{\bar w}}^*\right\|_2$. Note that this behavior is also typical of discrete entropy measures, as we stress in the next sections. Indeed, if the weights are more ``diverse'' to each other,  the distance w.r.t. the discrete uniform pmf $\mathbf{{\bar w}}^*$ is higher, the ESS and the entropy of ${\bar w}$ are smaller.
On the other hand, if the normalized weights are more similar to each other, they are all closer to the value $1/N$, so that the distance w.r.t. the discrete uniform pmf $\mathbf{{\bar w}}^*$ is smaller. As a consequence, the corresponding ESS and the entropy of ${\bar w}$ would be greater.
Hence, it appears natural to consider the possibility of using other discrepancy and/or entropy measures to design alternative ESS expressions. Highlighting these types of connections is relevant since (a) we can extend the range of applications of the ESS formulas (applying that expressions in other fields) and (b)  derivations employed in other fields can be used to design novel ESS formulas.
\newline
\newline
{\bf Why discrepancy measures.} The maximum ESS value is obtained when $\mathbf{{\bar w}}=\mathbf{{\bar w}}^*$, i.e., all the normalized weights are equal to $1/N$,
$$
\bar{w}_1=...=\bar{w}_N=\frac{1}{N}.
$$
This can be considered a good scenario (and confused with the optimal one)  since the ideal Monte Carlo estimator in Eq. \eqref{IdealMCestimator} can be interpreted as an  estimator with ``equal weighted samples'' (each $h({\bf x}_n)$ is multiplied by a factor $1/N$).
The problem is that with an IS scheme, the case $\bar{w}_1=...=\bar{w}_N=\frac{1}{N}$ (or $\bar{w}_1\approx...\approx\bar{w}_N\approx \frac{1}{N}$) can occur also in catastrophic scenarios, for instance, when all the samples are located in a tail of the target distribution (that often is a quite flat region), or when the samples are very close to each other. 
\newline
However, the ESS formulas above based on the discrepancy approach are able to detect other critical scenarios. 
For instance,  the minimum ESS value is reached when just one weight concentrates all the probability mass ($\bar{w}_i=1$ and $\bar{w}_j=0$ for $i\neq j$), which is a situation  to be avoided within particle filtering and sequential Monte Carlo schemes \cite{Aralumpalam02,Chopin02,Djuric03,Doucet00}. Therefore,
 this idea of this discrepancy approach has gained strength in the literature, and the ESS approximations above have been widely applied.     
\newline
\newline
In the following, we describe five conditions that a generic ESS approximation based only on the information of the normalized weights must satisfy. Then we show that the family of functions proposed in \citep{Huggins15} fulfills these five conditions. Furthermore, we link this G-ESS family  with the R\'enyi entropy providing also some theoretical results.

%%%%%%%%%%%%%%%%%%%%%%%%%%%%%%%%%%%%%%%%%%%%
%\section{Problem statement and background}
%%%%%%%%%%%%%%%%%%%%%%%%%%%%%%%%%%%%%%%%%%%%

%%%%%%%%%%%%%%%%%%%%%%%%%%%%%%%%%%%%%%%%%%%%
\section{Generalized ESS functions}\label{GESSsect}
%%%%%%%%%%%%%%%%%%%%%%%%%%%%%%%%%%%%%%%%%%%%
Considering the practical approach employed above for defining ESS formulas as discrepancy-diversity measures, here we describe the five properties that a generalized ESS measure (G-ESS) should satisfy, based only on the information of the normalized weights. We list below five conditions. The formulas which satisfy all of them can be applied as suitable ESS measures in practical applications (within IS or sequential IS schemes). Otherwise, if they satisfy at least the first three conditions, they can be considered as discrepancy measures with respect to the the uniform pmf, but have not the ability to be "particle/sample counters", as we clarify below with practical examples.
\newline
First of all, note that any possible G-ESS is a function of the vector of normalized weights ${\bf \bar{w}}=[\bar{w}_1,\ldots,\bar{w}_N]$,  
\begin{eqnarray}
\mbox{ESS}_N({\bf \bar{w}})= \mbox{ESS}_N(\bar{w}_1,\ldots,\bar{w}_N):  \mathcal{S}_{N}\rightarrow [1,N],
\end{eqnarray}
where $\mathcal{S}_{N} \subset \mathbb{R}^{N}$ represents the {\it unit simplex} in $\mathbb{R}^{N}$. Namely, the variables $\bar{w}_1,\ldots,\bar{w}_N$ are subjected to the following constraint:
 \begin{eqnarray}
 \bar{w}_1+\bar{w}_2+\ldots+\bar{w}_N=1.
 \end{eqnarray}
Moreover, we denoted
\begin{eqnarray}
 \label{VectMax}
{\bf \bar{w}}^*=\left[\frac{1}{N},\ldots,\frac{1}{N}\right],
\end{eqnarray}
and the vertices of the simplex $\mathcal{S}_{N}$ are denoted as
\begin{eqnarray}
 \label{VerticesSimplex}
 {\bf \bar{w}}^{(j)}&=&[\bar{w}_1=0,\ldots,\bar{w}_j=1,\ldots,\bar{w}_{N}=0], 
\end{eqnarray}
i.e., $\bar{w}_j=1$ and $\bar{w}_n=0$ (it can occur only if $\pi({\bf x}_n)=0$), for $n\neq j$ with $j\in\{1,\ldots,N\}$.  
\newline
\newline 
Below we list the five conditions that $\mbox{ESS}_N({\bf \bar{w}})$ should fulfill:
\begin{enumerate}
\item[C1.]  {\bf Symmetry:} $\mbox{ESS}_N$ must be invariant under any permutation of the weights, i.e.,
 \begin{equation}
 \label{SymmCond}
 \mbox{ESS}_N(\bar{w}_1,\bar{w}_2,\ldots,\bar{w}_N)= \mbox{ESS}_N(\bar{w}_{j_1},\bar{w}_{j_2},\ldots,\bar{w}_{j_N}),
\end{equation}
for any possible set of indices  $\{j_1,\ldots,j_N\}=\{1,\ldots,N\}$.
\item[C2.]  {\bf Maximum condition:} A maximum value is  $N$ and it is reached at ${\bf \bar{w}}^*$ (see Eq. \eqref{VectMax}), i.e.,
 \begin{equation}
 \mbox{ESS}_N\left({\bf \bar{w}}^*\right)=N\geq  \mbox{ESS}_N({\bf \bar{w}}).
 \end{equation}
%imply that $\mbox{ESS}_N$ has not local maxima  
 \item[C3.]  {\bf Minimum condition}: the minimum value is $1$ and it is reached (at least) at the vertices ${\bf \bar{w}}^{(j)}$ of the unit simplex  in Eq. \eqref{VerticesSimplex},
  \begin{equation}
 \mbox{ESS}_N({\bf \bar{w}}^{(j)})=1\leq  \mbox{ESS}_N({\bf \bar{w}}).
 \end{equation} 
 for all $j\in \{1,\ldots,N\}$.
  \item[C4.]  {\bf Unicity of extreme values:}  The maximum at ${\bf \bar{w}}^*$ is unique and the the minimum value $1$ is reached {\it only} at the vertices ${\bf \bar{w}}^{(j)}$, for all $j\in \{1,\ldots,N\}$.
\item[C5.] {\bf Stability of the rate} $\mbox{ESS}_N/N${\bf:} Consider the vector of weights ${\bf \bar{w}}\in \mathbb{R}^N$ and the vector 
${\bf \bar{v}}=[\bar{v}_1,\ldots,\bar{v}_{MN}] \in \mathbb{R}^{MN}, \quad  M\geq 1$,
obtained repeating and scaling by $\frac{1}{M}$ the entries of  ${\bf \bar{w}}$, i.e., 
\begin{equation}
\label{RepVect}
{\bf \bar{v}}=\frac{1}{M}[\underbrace{{\bf \bar{w}},{\bf \bar{w}},\ldots,{\bf \bar{w}}}_{M-times}].
\end{equation}
The invariance condition is expressed as
\begin{eqnarray}
\label{StableCon}
%\frac{\mbox{ESS}_N({\bf \bar{w}})}{N}&=&\frac{\mbox{ESS}_{MN}({\bf \bar{v}})}{MN} \nonumber \\
\mbox{ESS}_N({\bf \bar{w}})&=&\frac{1}{M} \mbox{ESS}_{MN}({\bf \bar{v}}),
\end{eqnarray}
%\frac{\mbox{ESS}_{mN}({\bar \alpha}_1=\frac{1}{m}\bar{w}_1,\ldots,{\bar \alpha}_N=\frac{1}{m}\bar{w}_N,{\bar \alpha}_{N+1}=\frac{1}{m}\bar{w}_1,\ldots,{\bar \alpha}_{mN}=\frac{1}{m}\bar{w}_N)}{mN}
for all $M \in \mathbb{N}^+$. 
%$$
%\mbox{ESS}_N({\bf \bar{w}})=\frac{1}{m} \mbox{ESS}_{mN}({\bf \bar{v}}).
%$$
\end{enumerate}
This last requirement can be interpreted as an adjustment of the  well-known {\it homogeneity} (scale-invariance) condition for real functions.\footnote{A function $f({\bf x})$ is said to be homogeneous of degree $k$ if $f(c \ {\bf x})=c^kf({\bf x})$ where $c$ is a non-zero constant value.} Note that, given conditions C2 and C3, we always have
\begin{equation}
1 \leq \mbox{ESS}_N({\bf {\bar w}})\leq N.
\end{equation}
If at least C1, C2 and C3 are fulfilled, the G-ESS 
can be considered a discrepancy measure with respect to the the uniform pmf. If also C4 is satisfied, then it is a proper discrepancy measure  since it reaches the maximum value (that is $N$) only at ${\bar w}^*$ and the minimum value (that is 1) only at the vertices ${\bf \bar{w}}^{(j)}$. However, if C5 is not ensured, the formula cannot be considered a useful ESS function from a practical point of view.
\newline
\newline
{\bf On the condition C5.} For clarifying this condition, consider the vector $\mathbf{\bar{v}} =[0,1,0]$ with $N=3$, and the two additional vectors obtained repeating $\mathbf{\bar{v}}$ two or three times,
$$
\begin{aligned}
%\mathbf{\bar{v}} & =[0,1,0], \\
\mathbf{\bar{v}}^{\prime} & =\left[0, \frac{1}{2}, 0,0, \frac{1}{2}, 0\right]=\frac{1}{2}[\mathbf{\bar{v}}, \mathbf{\bar{v}}], \\
\mathbf{\bar{v}}^{\prime \prime} & =\left[0, \frac{1}{3}, 0,0, \frac{1}{3}, 0,0, \frac{1}{3}, 0\right]=\frac{1}{3}[\mathbf{\bar{v}}, \mathbf{\bar{v}}, \mathbf{\bar{v}}],
\end{aligned}
$$
 We would like to obtain $\mbox{ESS}_3(\mathbf{\bar{v}})=1, \mbox{ESS}_{6}\left(\mathbf{\bar{v}}^{\prime}\right)=2$ and $\mbox{ESS}_{9}\left(\mathbf{\bar{v}}^{\prime \prime}\right)=3$, i.e., the ratio $\frac{\mbox{ESS}_N}{N}$ is constant, i.e.,
$$
\frac{\mbox{ESS}_3(\mathbf{\bar{v}})}{3}=\frac{\mbox{ESS}_{6}\left(\mathbf{\bar{v}}^{\prime}\right)}{6}=\frac{\mbox{ESS}_{9}\left(\mathbf{\bar{v}}^{\prime \prime}\right)}{9}=\frac{1}{3} .
$$
A more intuitive explanation is as follows. If we have a vector of normalized weights $\mathbf{\bar{v}}^{\prime} =\left[0, \frac{1}{2}, 0,0, \frac{1}{2}, 0\right]$, we would like to get $\mbox{ESS}_6(\mathbf{\bar{v}}')=2$, since we have $2$ effective samples instead of $6$ (at most we have $2$ effective samples; we can just say that, looking the vector $\mathbf{\bar{v}}^{\prime}$). Now, if we have a vector
$\mathbf{\bar{v}}^{\prime \prime}=\left[0, \frac{1}{3}, 0,0, \frac{1}{3}, 0,0, \frac{1}{3}, 0\right]$, we would like to obtain  $\mbox{ESS}_9(\mathbf{\bar{v}}^{\prime \prime})=3$. 
\newline
From another point of view, since $\mathbf{\bar{v}}^{\prime \prime}$ can be seen as $\mathbf{\bar{v}}^{\prime \prime}=\frac{1}{3}[\mathbf{\bar{v}},\mathbf{\bar{v}},\mathbf{\bar{v}}]$, where $\mathbf{\bar{v}} =[0,1,0]$, we would like that the ESS$_N$ formula would be able to count {\it effective samples} in the same way in different pieces of a vector. Namely, if $ESS_3(\mathbf{\bar{v}})=1$ and $\mathbf{\bar{v}}^{\prime \prime}$ is formed by three repetitions of $\mathbf{\bar{v}}$  then we expect to obtain $\mbox{ESS}_9(\mathbf{\bar{v}}^{\prime \prime})=3$.
 Any result that differs from these ones, does not make sense from a practical point of view, e.g., within a particle filter or sequential Monte Carlo scheme. 
%{\color{magenta} ability of detection...  (non-zero elements, no se)... and stability as $N$ grows... keeping the same cataracterstics/features if we repeat the same vector in order to build a larger one... }
\newline
{\bf Classification of G-ESS.} Given the previous observations, we can provide a classification of the possible G-ESS functions. Table \ref{TableClass} classifies the G-ESS functions in different families depending on the conditions fulfilled. Table \ref{TableClass} shows the cases found in different families of ESS measures \cite{ESSarxiv16}.  Recall that the first three conditions are strictly required, to be considered a discrepancy measure with respect to the uniform pmf. All the G-ESS functions which satisfy at least the first four conditions, i.e., from C1 to C4, are called {\it proper} functions. If all the conditions are fulfilled they are called {\it proper and }stable.
We are interested in this last type of G-ESS expressions, proper and stable.
%\newline
%\newline
{\rem Only the proper and stable G-ESS functions are  useful from a practical point of view, in order to be employed as ESS measures.}

\begin{table}[h!] 
\caption{Classification of G-ESS formulas.}\label{TableClass}
\begin{center}
\begin{tabular}{|c|c|c|c|c|c|}
\hline Class of G-ESS & C1 & C2 & C3 & C4 & C5 \\
\hline \hline Degenerate  & \checkmark & \checkmark & \checkmark & x &x \\
\hline Proper  & \checkmark & \checkmark & \checkmark & \checkmark & x \\
\hline Degenerate and Stable  & \checkmark & \checkmark & \checkmark & x & \checkmark \\
\hline Proper and Stable  & \checkmark & \checkmark & \checkmark & \checkmark & \checkmark \\
\hline
\end{tabular}
\end{center}
\end{table}

\noindent
In order to clarify the previous remark, and the importance of the five conditions, below we show the importance of the condition C5, introducing a family that fulfills the first 4 conditions but does not satisfy the last condition.
%%%%%%%%%%%%%%%%%%%%%%%%%%%%%%%%%%%%%%%%%%%%
\paragraph{Example of a proper but non-stable G-ESS family.}
%%%%%%%%%%%%%%%%%%%%%%%%%%%%%%%%%%%%%%%%%%%%
%\newline
%\newline
Here, as an example, we introduce a G-ESS family such that the formulas in that family are all proper but {\it not} stable. This means that all the contained G-ESS expressions can be used as discrepancy measure with respect to the uniform pmf but are not suitable to be employs as ESS measures (within particle filters or sequential Monte Carlo schemes).
We can design a G-ESS family based on the $L_p$ distance between ${\bf \bar{w}}$ and ${\bf \bar{w}}^*$ which satisfies the first four conditions above. This could be an intuitive idea when we are interested in discrepancy measures. We can in fact define the family
\begin{align}\label{BasedOndistance}
\mbox{ESS-D}_N^{(p)}({\bf \bar{w}})=\frac{1}{ \alpha_p ||{\bf \bar{w}}-{\bf \bar{w}}^* ||_p +\frac{1}{N}}, \quad \quad \alpha_p=\frac{N-1}{N\left[\frac{N-1+(N-1)^p}{ N^p}\right]^{1 / p}},
\end{align}
where
$$
||{\bf \bar{w}}-{\bf \bar{w}}^* ||_p=\left(\sum_{i=1}^N \left|{\bar w}_n-\frac{1}{N}\right|^p\right)^{1/p}, \qquad p>0.
$$
It is possible to show that  $\mbox{ESS-D}_N^{(p)}({\bf \bar{w}}^*)=N$ and
$\mbox{ESS-D}_N^{(p)}({\bf \bar{w}}^{(j)})=1$. Hence, $\mbox{ESS-D}^{(p)}$ fulfills C1, C2, and C3 and it is also easy to show that satisfies C4. 
Hence, this family can be employed as a discrepancy measure with respect to the uniform pmf ${\bf \bar{w}}^*$. However, it is not a good ESS measure since  does not satisfy C5 (it is not stable). For clarifying this point, let us consider some examples comparing $\mbox{ESS-D}_N^{(p)}$ with $p=2$ in  Eq. \eqref{BasedOndistance},  with other ESS formulas (but {\it stable}) in  Eqs. \eqref{FirstDef_P} and \eqref{MaxEq}. The results are given in Table \ref{TableESSex}.

\begin{table}[h!]
\footnotesize
\caption{Examples of ESS measures with different vector ${\bf \bar{w}}$ with dimension $N=5$. Note that $\mbox{ESS-D}_N^{(p)}({\bf \bar{w}})$ with $p=2$ in  Eq. \eqref{BasedOndistance} is proper but non-stable, whereas the rest of  two ESS formulas in  Eqs. \eqref{FirstDef_P}-\eqref{MaxEq} are both proper and stable.}
\label{TableESSex} 
\vspace{-0.2cm}
\begin{center}
\begin{tabular}{|c|c|c|c|c|c|}
%\hline
%{\bf Parameter:} & ${\bf  r\rightarrow 0}$ & ${\bf  r\rightarrow 1}$&  ${\bf  r\rightarrow \infty}$ \\
\hline
& \multirow{2}{*}{(a)} & \multirow{2}{*}{(b)} & \multirow{2}{*}{(c)} & \multirow{2}{*}{(d)} & \multirow{2}{*}{(e)} \\
  & & & & & \\
  ${\bf \bar{w}} \Longrightarrow$    & $[1,0,0,0,0]$ & $\left[\frac{1}{2},\frac{1}{2},0,0,0\right]$  & $\left[\frac{1}{3},\frac{1}{3},\frac{1}{3},0,0\right]$ Ê& $\left[\frac{1}{4},\frac{1}{4},\frac{1}{4},\frac{1}{4},0\right]$ &  $\left[\frac{1}{5},\frac{1}{5},\frac{1}{5},\frac{1}{5},\frac{1}{5}\right]=\bar{{\bf w}}^*$\\
  & & & & & \\
\hline
  & & & & & \\
$\mbox{ESS-D}_5^{(2)}({\bf \bar{w}})$  ---  Eq. \eqref{BasedOndistance} & 1 &  1.45 &  1.90 &  2.5  & 5 \\
  & & & & & \\
\hline
  & & & & & \\
$\frac{1}{\sum_{n=1}^5 {\bar w}_n^2}$  ---  Eq. \eqref{FirstDef_P}   & 1 &  2 & 3 & 4 & 5 \\
  & & & & & \\
\hline
  & & & & & \\
$\frac{1}{\max \bar{w}_n}$  ---  Eq. \eqref{MaxEq}   & 1 &  2 & 3 & 4 & 5 \\
  & & & & & \\
\hline
\end{tabular}
\end{center}
\end{table} 

\noindent
the formula $\mbox{ESS-D}_N^{(p)}$ does not provide the desired results, with the exception of the cases (a) and (e) (the first and the last scenarios) that are related to the condition C2 and C3. Namely, $\mbox{ESS-D}_N^{(p)}$ is not a good particle counter, unlike the other two ESS formulas. For instance, in case of ${\bf \bar{w}}=[\bar{w}_1=\frac{1}{3},\bar{w}_2=\frac{1}{3},\bar{w}_3=\frac{1}{3},\bar{w}_4=0,\bar{w}_5=0]$ using  {\it just}  the information of these normalized weights $\bar{w}_n$, we can just assert that we have three effective samples, whereas  $\mbox{ESS-D}_5^{(2)}({\bf \bar{w}})$  returns $\approx 1.90$.

%%%%%%%%%%%%%%%%%%%%%%%%%%%%%%%%%%%%%%%%%%%%
\section{Huggins-Roy's ESS family}\label{HR_Sect}
%%%%%%%%%%%%%%%%%%%%%%%%%%%%%%%%%%%%%%%%%%%%
The Huggins-Roy's ESS family introduced in \citep{Huggins15} is defined as
\begin{eqnarray}
\mbox{ESS-H}_N^{(\beta)}({\bf \bar{w}})&=&\left(\frac{1}{\sum_{n=1}^N {\bar w}_n^\beta}\right)^{\frac{1}{\beta-1}}, Ê\\
&=&\left(\sum_{n=1}^N {\bar w}_n^\beta \right)^{\frac{1}{1-\beta}}, Ê\quad \quad \beta\geq 0.
\end{eqnarray}
Table \ref{TableStableGESS}  shows below that the Huggins-Roy's family contains all the most important, {\it proper and stable} G-ESS functions introduced in literature.
The special cases with $\beta=0$ and $\beta=1$  bring to two undetermined expressions that will be solved and clarified below (when the relationship with R\'enyi entropy is shown). 
We can easily note that $1\leq \mbox{ESS-H}_N^{(\beta)}({\bf \bar{w}})\leq N$ for all $\beta\geq 0$. More generally, it is possible to observe that for $\beta\neq 0$ the conditions C1, C2, C3 and C4 are fulfilled (with the exception of $\beta=0$ that does not satisfy C4). Furthermore, the condition C5 is also satisfied, for all $\beta$, as we show next.
\newline
\newline
{\bf Proof.} In order to prove that  C5 is  satisfied, for simplicity let us consider a vector ${\bf \bar{v}}=\frac{1}{2}[{\bf \bar{w}},{\bf \bar{w}}]$, defined repeating twice the vector ${\bf \bar{w}}$ (i.e., $M=2$). In this case, we have
\begin{eqnarray}
\mbox{ESS-H}_{2N}^{(\beta)}({\bf \bar{v}})&=&\left(\frac{1}{2^\beta}\sum_{n=1}^N {\bar w}_n^\beta+\frac{1}{2^\beta}\sum_{n=1}^N {\bar w}_n^\beta \right)^{\frac{1}{1-\beta}}, \nonumber  \\
&=& \left(\frac{1}{2^{\beta-1}}\sum_{n=1}^N {\bar w}_n^\beta \right)^{\frac{1}{1-\beta}}, \nonumber \\
&=&2\left(\sum_{n=1}^N {\bar w}_n^\beta \right)^{\frac{1}{1-\beta}}, \nonumber \\
&=& 2 \ \mbox{ESS-H}_{N}^{(\beta)}({\bf \bar{w}}), \quad \forall \beta,
\end{eqnarray}
which is exactly the condition in Eq. \eqref{StableCon}. The proof can be easily repeated for a value $M>2$.
\newline
\newline
{\rem Hence, all G-ESS functions (except for $\beta\rightarrow 0$) belonging to the Huggins-Roy's ESS family are {\it proper and stable}. For $\beta\rightarrow 0$, the corresponding ESS is degenerate and stable. 
Moreover, some specific cases provided in Table \ref{TableStableGESS}, coincide with other proper and stable  G-ESS  formulas proposed in \citep{ESSarxiv16}.}

\begin{table}[hbt]
\setlength{\tabcolsep}{4pt}
\caption{Relevant special cases contained in the Huggins-Roy's family. They are all proper and stable, except for $N-N_Z$ that is degenerate and stable.}
\label{TableStableGESS} 
\vspace{-0.2cm}
\begin{center}
\begin{tabular}{|c|c|c|c|c|}
%\hline
%{\bf Parameter:} & ${\bf  r\rightarrow 0}$ & ${\bf  r\rightarrow 1}$&  ${\bf  r\rightarrow \infty}$ \\
\hline
   $\beta\rightarrow 0$ &$\beta=1/2$ & $\beta \rightarrow 1$  & $\beta=2$Ê& $\beta  \rightarrow \infty$\\
\hline
\hline

\multirow{2}{*}{$N-N_Z$ }  & \multirow{2}{*}{$\left(\sum_{n=1}^N \sqrt{{\bar w}_n}\right)^2$} & \multirow{2}{*}{$\exp\left(-\sum_{n=}^N{\bar w}_n \log {\bar w}_n\right)$} & \multirow{2}{*}{$\frac{1}{\sum_{n=1}^N {\bar w}_n^2}$}  &   \multirow{2}{*}{$\frac{1}{\max[{\bar w}_1,\ldots,{\bar w}_N]}$ }    \\
& &  &  &\\
{\it where $N_Z$ is}&  & &   &   \\
{\it the number of}& & {\it Perplexity - Eq. \eqref{PerpEq}} &  {\it Standard formula}   &  {\it In Eq. \eqref{MaxEq}} \\
{\it zeros in ${\bf {\bar w}}$}& & \citep{pmc-cappe08,Robert10b}  & {\it  in Eq. \eqref{FirstDef_P}} - \cite{Kong92}  &   \cite{ESSarxiv16}  \\
%\hline
%\hline
%{\it sub-linear} & {\it super-linear} & {\it super-linear}  \\ 
\hline
%degenerate & proper   & & & \\
%and stable  & and stable    & and stable  & and stable  & and stable  \\
%\hline
\end{tabular}
\end{center}
\end{table} 

%%%%%%%%%%%%%%%%%%%%%%%%%%%%%%%%%%%%%%%%%%%%
\section{Relationship with the entropy measures}\label{EntropiesSect}
\subsection{Relationship with the R\'enyi entropy}
%%%%%%%%%%%%%%%%%%%%%%%%%%%%%%%%%%%%%%%%%%%%
In this section, we show the connection between the R\'enyi entropy and  Huggins-Roy's family.
The R\'enyi entropy \citep{Cover91} is defined as 
\begin{eqnarray}
R_N^{(\beta)}({\bf \bar{w}})&=&\frac{1}{1-\beta} \log \left[ \sum_{n=1}^N {\bar w}_n^\beta\right],  Ê\quad \quad \beta >0,Ê
\end{eqnarray}
Then, first noting that $\frac{1}{1-\beta} \log \left[ \sum_{n=1}^N {\bar w}_n^\beta\right]=\log \left[ \sum_{n=1}^N {\bar w}_n^\beta\right]^{\frac{1}{1-\beta} }$ and  taking the exponential of both sides of the equation above, we obtain  
\begin{equation}\label{ExpOfaEntropy}
\mbox{ESS-H}_N^{(\beta)}({\bf \bar{w}})=\exp\left(R_N^{(\beta)}({\bf \bar{w}})\right)=\left(\sum_{n=1}^N {\bar w}_n^\beta \right)^{\frac{1}{1-\beta}}, Ê\quad \quad \beta > 0.
\end{equation}
 In ecology, the exponential of the R\'enyi entropy defines the so-called {\it diversity indices} \citep{Jost06}. This means that   the Huggins-Roy's family contains and coincides with all the diversity indices derived by the R\'enyi entropy  \citep{Cover91,Jost06}. See Section \ref{EcoSect}, for further details. 
\newline
 Note that, for $\beta=0$, we have $R_N^{(0)}({\bf \bar{w}})=\log (N-N_Z)$  where $N_Z=\#\big\{\mbox{all ${\bar w}_n$:} \quad {\bar w}_n=0, \quad \forall n=1,\ldots,N \big\}$ (see \citep{Cover91} for further details), so that  $\mbox{ESS-H}_N^{(0)}({\bf \bar{w}})=N-N_Z$, as also shown in Table \ref{TableStableGESS}. For $\beta=1$, we have $R_N^{(0)}({\bf \bar{w}})=-\sum_{n=}^N{\bar w}_n \log {\bar w}_n$ \citep{Cover91}  then
\begin{equation}
\mbox{ESS-H}_N^{(1)}({\bf \bar{w}})=\exp\left(-\sum_{n=}^N{\bar w}_n \log {\bar w}_n\right),
\end{equation}
that is the perplexity in Eq. \eqref{PerpEq} \citep{pmc-cappe08,Robert10b}.  
%\newline
%\newline

\subsubsection{Inequalities for the G-ESS within Huggins-Roy family}

  One of advantages of the connection with the R\'enyi entropy is that we can obtain easily some theoretical results about $\mbox{ESS-H}_N^{(\beta)}$. Indeed, for instance, it is well-known that  \citep{Cover91}
$$
R_N^{(0)}({\bf \bar{w}})\geq R_N^{(1)}({\bf \bar{w}}) \geq R_N^{(2)}({\bf \bar{w}})\geq \ldots  R_N^{(\beta')}({\bf \bar{w}}) \ldots   \geq R_N^{(\infty)}({\bf \bar{w}}), \qquad \beta'\geq 2.
$$
 Then, since $\mbox{ESS-H}_N^{(\beta)}$ is an increasing monotonic function of $R_N^{(\beta)}$, we can also assert
\begin{equation}
\mbox{ESS-H}_N^{(0)}({\bf \bar{w}})\geq \mbox{ESS-H}_N^{(1)}({\bf \bar{w}}) \geq \mbox{ESS-H}_N^{(2)}({\bf \bar{w}})\geq \ldots  \mbox{ESS-H}_N^{(\beta')}({\bf \bar{w}}) \ldots    \geq \mbox{ESS-H}_N^{(\infty)}({\bf \bar{w}}).
\end{equation}
Namely, we can re-write
\begin{align}
\mbox{ESS-H}_N^{(\infty)}({\bf \bar{w}})&\leq \mbox{ESS-H}_N^{(\beta)}({\bf \bar{w}})  \leq \mbox{ESS-H}_N^{(0)}({\bf \bar{w}}), \nonumber  \\
 \frac{1}{\max \bar{w}_n}&\leq \mbox{ESS-H}_N^{(\beta)}({\bf \bar{w}})  \leq N-N_Z, \quad \beta\geq 0.
\end{align}
 Moreover, since from \citep{Cover91} we have
$$
R_N^{(2)}({\bf \bar{w}}) \leq 2 \mbox{ } ÊR_N^{(\infty)}({\bf \bar{w}}),
$$
we can also write 
\begin{equation}
\mbox{ESS-H}_N^{(2)}({\bf \bar{w}}) \leq 2 \mbox{ }  \mbox{ESS-H}_N^{(\infty)}({\bf \bar{w}}).
\end{equation}

%%%%%%%%%%%%%%%%%%%%%%%%%%%%%%%%%%%%%%%%%%%%
\subsection{Relationship with the Tsallis entropy}\label{TsallisSect}
%%%%%%%%%%%%%%%%%%%%%%%%%%%%%%%%%%%%%%%%%%%%
Another famous entropy family is the so-called Tsallis entropy \cite{Tsallis1988} (as known as $q$-logarithmic entropy \cite{Jost06}), defined as
\begin{align}
T_N^{(\alpha)}({\bf \bar{w}})&=   \frac{1}{\alpha-1} \left[1- \sum_{n=1}^N {\bar w}_n^\alpha\right],  Ê\quad \quad  \alpha > 0.
\end{align}
We can obtain a corresponding G-ESS family based on the Tsallis entropy, % taking the exponential of $T_N^{(\alpha)}({\bf \bar{w}})$ 
 after some additional simple operations of translation and scaling, i.e.,
\begin{align}
\mbox{ESS-T}_N^{(\alpha)}({\bf \bar{w}})&=\frac{(\alpha-1)(N-1)}{N^{1-\alpha}-1} T_N^{(\alpha)}({\bf \bar{w}})+1, \\
&=\frac{(\alpha-1)(N-1)}{N^{1-\alpha}-1}  
  \left[1- \sum_{n=1}^N {\bar w}_n^\alpha\right]+1,  Ê\quad \quad \alpha  > 0.
\end{align}
Note that
$$
1\leq \mbox{ESS-T}_N^{(\alpha)}({\bf \bar{w}}) \leq N.
$$
{\bf Special cases.} For $\alpha\rightarrow 0$, we get again the following degenerate and stable formula $\mbox{ESS-T}_N^{(0)}({\bf \bar{w}})=N-N_Z$, where $N_Z=\#\big\{\mbox{all ${\bar w}_n$:} \quad {\bar w}_n=0, \quad \forall n=1,\ldots,N \big\}$. For $\alpha\rightarrow \infty$, we have the degenerate expression $\mbox{ESS-T}_N^{(\infty)}({\bf \bar{w}})=N$ if ${\bf \bar{w}}\neq {\bf \bar{w}}^{(j)}$ for all $j \in\{1, \ldots, N\}$, or $\mbox{ESS-T}_N^{(\infty)}({\bf \bar{w}})=1$ if ${\bf \bar{w}}={\bf \bar{w}}^{(j)}$,  for all $j \in\{1, \ldots, N\}$.
\newline
Setting $\alpha=2$, we have
\begin{align}
\mbox{ESS-T}_N^{(2)}({\bf \bar{w}})%&=\frac{N(N-1)}{(N-1)}  
  %\left(1- \sum_{n=1}^N {\bar w}_n^2\right)+1,  \nonumber \\
 &=N 
  \left(1- \sum_{n=1}^N {\bar w}_n^2\right)+1,  \nonumber  \\
 &=N \ \mbox{Gini-impurity}({\bf{\bar w}})+1,   
\end{align}
where  we have used the definition of the function below,
\begin{align}
\mbox{Gini-impurity}({\bf{\bar w}})=1- \sum_{n=1}^N {\bar w}_n^2,
\end{align}
that is the so-called {\it Gini impurity} or {\it Gini's diversity index} or also known as {\it Gini-Simpson index} in biodiversity studies,  that is widely used in  machine learning within decision tree algorithms \cite{bishop2007,Krzywinski2017}. Moreover, from an ecology point of view, that $\mbox{Gini-impurity}({\bf{\bar w}})$  represents the probability that two individuals chosen at random are of different species. The Gini impurity is associated with the name of Edward H. Simpson, who introduced it as an index of diversity in 1949 \cite{SIMPSON1949}.
 Then, Corrado Gini used the formula (called as ``Gini impurity'') above in economics, statistics, and demography \cite{Gini12}. It is such a natural quantity that it has been used in many different fields and admits an unbiased estimator. Despite all these benefits, $\mbox{Gini-impurity}({\bf{\bar w}})$ is not directly an effective number and needs an additional translation and scaling, becoming $\mbox{ESS-T}_N^{(2)}({\bf \bar{w}})$. Moreover, the final expression is not stable.
\newline
\newline
It is also interesting to remark that the final form of $\mbox{ESS-T}_N^{(\alpha)}({\bf \bar{w}})$ resembles the G-ESS family $\mbox{ESS-V}_N^{(r)}(\mathbf{\bar{w}})$ introduced in \cite{ESSarxiv16}, 
$$
\mbox{ESS-V}_N^{(r)}(\mathbf{\bar{w}})=\frac{N^{r-1}(N-1)}{1-N^{r-1}}\left[\sum_{n=1}^N \bar{w}_n^r\right]+\frac{N^r-1}{N^{r-1}-1},
\qquad r>0,
$$
%that contains $\mbox{ESS-S}_N^{(1/2)}(\mathbf{\bar{w}})=\left(\sum_{n=1}^N \sqrt{{\bar w}_n}\right)^2$ for $r=1/2$, that is a proper and stable ESS formula. 
However, generally the ESS expressions contained in $\mbox{ESS-V}_N^{(r)}(\mathbf{\bar{w}})$ and $\mbox{ESS-T}_N^{(\alpha)}(\mathbf{\bar{w}})$ are not stable. 
For this reason, in this work we focus on  mainly Huggins-Roy ESS family.
Furthermore, it is also possible to find another transformation, instead of the standard exponential function $\exp(\cdot)$ (as for the R{\'e}nyi entropy), that converts the Tsallis entropy into the  Huggins-Roy ESS family. That is the so-called $q$-exponential function \cite{Cover91}:
\begin{align}\label{ExtendedExp}
\exp_\alpha(t)= \begin{cases}(1+(1-\alpha) t)^{1 /(1-\alpha)} & \text { if } \alpha \neq 1, \\ \exp(t) & \text { if } \alpha=1.\end{cases}
\end{align}
After some manipulations, we arrive to
\begin{align}
\exp_\alpha\left(T_N^{(\alpha)}({\bf \bar{w}})\right)=\mbox{ESS-H}_N^{(\alpha)}(\mathbf{{\bar w}}).
\end{align}
Hence, this confirms in a generalized-sense the definition  of a diversity index as ``exponential of an entropy'' given in Eq. \eqref{ExpOfaEntropy}, used in ecology.

%%%%%%%%%%%%%%%%%%%%%%%%%%%%%%%%%%%%%%%%%%%%
\section{Other stable G-ESS expressions}\label{other_ESS_Stable}
%%%%%%%%%%%%%%%%%%%%%%%%%%%%%%%%%%%%%%%%%%%%
All the ESS formulas contained in the Huggins-Roy family are proper and stable, as we have shown in Section \ref{HR_Sect}. The converse statement is not true, i.e., there are other proper and stable G-ESS formulas that are not contained in the Huggins-Roy family. We provide some examples below. 
\newline
\newline
{\bf Another degenerate and stable formula.} We start with an additional example of degenerate and stable expression:
\begin{align}\label{Nplus}
\mbox{ESS-Plus}_{N}({\bar {\bf w}})=N^{+}=\#\left\{\bar{w}_n \geq 1 / N, \quad \forall n=1, \ldots, N\right\}.
\end{align}
It represents the number of the normalized weights bigger or equal to $1/N$. This ESS expression is stable but degenerate. The issue is that $\mbox{ESS-Plus}_{N}({\bar {\bf w}})$ reaches the minimum value 1 even at points that are not the vertices ${\bar {\bf w}}^{(j)}$ of the simplex (see Eq. \eqref{VerticesSimplex}). For instance, with ${\bar {\bf w}}=[0.8,0,0.2]$ we get $\mbox{ESS-Plus}_{3}({\bar {\bf w}})=1$, but we would like to reach the minimum value 1 only at the vertex  ${\bar {\bf w}}^{(1)}=[1,0,0]$, ${\bar {\bf w}}^{(2)}=[0,1,0]$ and  ${\bar {\bf w}}^{(3)}=[0,0,1]$. However, $\mbox{ESS-Plus}_{N}({\bar {\bf w}})$ is much more useful than another degenerate and stable formula that we already found in Table \ref{TableStableGESS}, i.e., $N-N_Z$. Indeed,  $N-N_Z$ is degenerate since reaches the maximum value, $N$, in any ${\bar {\bf w}}$ that does not contain any zero (instead of only at ${\bar {\bf w}}^*$). This makes $N-N_Z$ much less useful from a practical point of view, for instance, within a particle filter. Whereas  
$\mbox{ESS-Plus}_{N}({\bar {\bf w}})$ could be perfectly employed within a particle filter, considering it as a more conservative ESS formula with respect to other ESS expressions.
\newline
\newline
{\bf Other proper and stable formulas.} Let us define
\begin{align}
\left\{\bar{w}_1^{+}, \ldots, \bar{w}_{N^{+}}^{+}\right\}=\left\{\mbox{all } \bar{w}_n \mbox{ such that: }  \bar{w}_n \geq 1 / N, \quad \forall n=1, \ldots, N\right\},
\end{align}
where $N^+$ is given in Eq. \eqref{Nplus}, i.e., $N^{+}=\# \left\{\bar{w}_1^{+}, \ldots, \bar{w}_{N^{+}}^{+}\right\}$. Now, it is possible to define a correct-proper version of $\mbox{ESS-Plus}$  \cite{ESSarxiv16}, i.e., 
\begin{align}
\mbox{ESS-Q}_N(\mathbf{{\bar w}})&=-N \sum_{i=1}^{N^{+}} \bar{w}_i^{+}+N^{+}+N,  \nonumber \\
&=N^{+}+N\left(1-\sum_{i=1}^{N^{+}} \bar{w}_i^{+}\right), \nonumber \\    
&=N^{+}+N\left(\sum_{i=1}^{N-N^{+}} \bar{w}_i^{-}\right)=N^{+}+N\gamma,\label{Qplus}
\end{align}
where $\bar{w}_i^{-}$ are all the normalized weights such that $<1/N$, and $\gamma=\sum_{i=1}^{N-N^{+}} \bar{w}_i^{-} \leq  1$. Note that $\gamma=0$ in the two extreme cases ${\bar {\bf w}}={\bar {\bf w}}^{(j)}$ and ${\bar {\bf w}}={\bar {\bf w}}^*$ and $\mbox{ESS-Q}_N(\mathbf{{\bar w}})=N^{+} +0=N^{+}$, i.e.,  we have $\mbox{ESS-Q}_N(\mathbf{{\bar w}}^{j})=1$ and $\mbox{ESS-Q}_N(\mathbf{{\bar w}}^*)=N$ as expected.
In all the other scenarios, a portion of all the number of samples (that is $\gamma N$ with $\gamma \leq  1$) is added to $N^{+}$. The resulting ESS formula is proper and stable. This measure is also related to $L_1$ between the two pmfs \cite{ESSarxiv16}.
\newline
\newline
Another proper and stable ESS expression introduced in the literature  is based on the Gini inequality coefficient, widely applied in economics \cite{Gini21,ESSarxiv16}.  First of all, we define the non-decreasing sequence of normalized weights as
\begin{align}
\bar{w}_{(1)} \leq \bar{w}_{(2)} \leq \ldots \leq \bar{w}_{(N)},
\end{align}
obtained sorting in ascending order the entries of the vector ${\bar {\bf w}}$.  Let us consider the Gini inequality coefficient $G({\bar {\bf w}})$ introduced in economy for measuring the wealth inequality  can be defined as follows \cite{Gini1921,Gini12,Gini21},
\begin{align}
G({\bar {\bf w}})=2 \frac{s({\bar {\bf w}})}{N}-\frac{N+1}{N}, \quad \mbox{ where } \quad s({\bar {\bf w}})=\sum_{n=1}^N n \bar{w}_{(n)}.
\end{align}
It is not the unique formulation: there are various equivalent formulations of the Gini coefficient \cite{Gini13,Lorenz05}. Then, the corresponding G-ESS function is given by
\begin{align}
\mbox{ESS-Gini}_N({\bar {\bf w}})&=-N G({\bar {\bf w}})+N, \nonumber \\
&=-N \left[2 \frac{s({\bar {\bf w}})}{N}-\frac{N+1}{N}\right]+N,  \nonumber\\
&=-2s({\bar {\bf w}})+1+2N, \nonumber \\
&=-2\sum_{n=1}^N n \bar{w}_{(n)}+1+2N, \label{GINIess}
\end{align}
which is proper and stable. It can be easily also shown that $\mbox{ESS-Gini}_N({\bar {\bf w}}^{(j)})=-2N+1+2N=1$ for all $j$ and  $ \mbox{ESS-Gini}_N({\bar {\bf w}}^*)=-2\frac{1}{N}\frac{N(N+1)}{2}+1+2N=N$.
\newline
\newline
 The fact that some proper and stable ESS expressions do not belong to the Huggins-Roy family shows there is still room and need for further research on in this topic. For instance, new proper and stable formulas could be discovered, as we shown in the next section.

\section{Connections with other research fields: extended range of applications}\label{EcoSect0}

In the previous section, we have already seen that  the connections with the R\'enyi and Tsallis entropy families  show  the existence  of other relationships with many studies in different fields (e.g., ecology and machine learning). The benefit  of creating these bridges between fields is bidirectional: different ideas used in other fields can be applied as ESS in a IS context and, vice-versa, ESS formulas proposed for IS could be employed in other fields. 
\newline
A clear example of this benefit  is given, in Section \ref{ESSpolSect}, where we discover a new proper and stable ESS formula, that has been introduced in political science. %This section describes more in details this extended range of applications.

\subsection{ESS in ecology}\label{EcoSect}

The connection with the R\'enyi entropy shows that the G-ESS functions of the Huggins-Roy's family  are also diversity indices \citep{Jost06}.  More specifically, the exponential of the R\'enyi entropy is known in ecology as the {\it Hill number} of order $\beta$ \citep{Jost06}.
The Hill numbers are the most important measures of biological diversity.
For instance, the Hill number of order 0 corresponds to $\mbox{ESS-H}_N^{(0)}({\bf \bar{w}})=N-N_Z$, and represents the number of species. This is also called the {\it species richness} in ecology, and is often used as a measure of diversity in the popular media and the ecology literature. However, it does not make any distinction between a rare species and a common species. Moreover, $\mbox{ESS-H}_N^{(0)}$ does not provide any information about the balance between the species that are involved.
\newline
In Section \ref{TsallisSect}, we have seen that the formula  $1- \sum_{n=1}^N {\bar w}_n^2$ is  called Gini-impurity in machine learning. Whereas,  in ecology, it  is called  Gini-Simpson index, since Simpson introduced it as an index of diversity \cite{SIMPSON1949}.  Moreover, since the sum of the squares $\sum_{n=1}^N {\bar w}_n^2$ can be interpreted as a {\it measure of concentration}  (see Section \ref{ConcMeSect}), the Hill number (diversity) of order 2, $\mbox{ESS-H}_N^{(2)}(\mathbf{{\bar w}})$, is also called the {\it inverse Simpson concentration} in ecology  \cite{SIMPSON1949}. Furthermore, 
the diversity of order $\infty$, i.e.,
$\mbox{ESS-H}_N^{(\infty)}(\mathbf{{\bar w}})=1 / \max {\bar w}_i$,
is known as {\it the Berger-Parker index} in ecology \cite{Berger70}. 
While the Hill number of order 0 gives rare species the same importance as any other, diversity of order $\infty$ ignores them and takes into account only the dominant species. More generally, the parameter $\beta$ controls the sensitivity of the diversity measure $\mbox{ESS-H}_N^{(\beta)}$ to rare species, with higher values of $\beta$ corresponding to measures less sensitive to rare species. In other words, $\beta$  reflects the inverse of the importance given to rare species.

%%%%%%%%%%%%%%%%%
\subsection{ESS in economics}\label{ConcMeSect}
%%%%%%%%%%%%%%%%%

The ESS indices have been also widely employed (under other names) as metrics for portfolio dispersion and/or concentration. The effective number of positions held in a portfolio is usually measured as  $\mbox{ESS}_N({\bf {\bar w}})=\frac{1}{\sum_{n=1}^N {\bar w}_n^2 }$, where the normalized weights ${\bar w}_n$ represent the proportion of market value invested in each security. 
A high value of $\mbox{ESS}_N({\bf {\bar w}})$  implies a very diversified portfolio (at most different $N$ equally weighted positions).
The formula $\frac{1}{\sum_{n=1}^N {\bar w}_n^2}$ has been shown to be one of the most efficient measures of portfolio diversification.  It has been also used as a constraint to force a portfolio to hold a minimum number of effective assets, denoted for instance as $N_{\text {eff }}$  (e.g., $\|{\bf {\bar w}}\|^2 \leq N_{\text {eff }}^{-1}$).
\newline
{\bf Concentration measures.} In Section \ref{EcoSect}, we have seen that there is a coincidence between Hill numbers and the Huggins-Roy ESS formulas. More generally, the reciprocals of the Hill numbers (hence, the reciprocals of the  G-ESS formulas as well) have been used in economics as concentration measures, i.e., 
\begin{equation}
\mbox{Conc}_N^{(\beta)}({\bf \bar{w}})=\frac{1}{ \mbox{ESS-H}_N^{(\beta)}({\bf \bar{w}})}=\left(\sum_{n=1}^N {\bar w}_n^\beta \right)^{\frac{1}{\beta-1}}, Ê\quad \quad \beta > 0.
\end{equation}
As an example, we could investigate if  an industry or a market is concentrated in the hands of a small number of large players. Let assume there are $N$ competing companies in a given industry, each one occupying a portion of the market represented by the normalized weights $\bar{w}_1, \ldots, \bar{w}_N$, then the concentration $1 / \mbox{ESS-H}_N^{(\beta)}({\bf \bar{w}})$ is maximized when one company has a monopoly, i.e., when ${\bf \bar{w}}={\bf \bar{w}}^{(j)}$ (the $j$-th company has conquered all the market, $\bar{w}_j=1$). Namely, a concentration index  ranges from $1/N$ (in case of perfect competition) to $1$ (in case of monopoly), where $N$ represents the number of companies in the market.
The concentration measure for $\beta=2$, i.e., $\mbox{Conc}_N^{(2)}({\bf \bar{w}})=\frac{1}{ \mbox{ESS-H}_N^{(2)}({\bf \bar{w}})}=\sum_{n=1}^N {\bar w}_n^2 $ is known as {\it Herfindahl-Hirschman index} in economy. Finally, in the previous section, we have also seen the application of similar indices for measuring the wealth inequality, e.g., using  the Gini coefficient \cite{Gini12,Gini1921,Gini13}. 
%{\color{magenta}
%See Hannah and Kay [134] or Chakravarty and Eichhorn [63], for instance.

%\url{https://en.wikipedia.org/wiki/HerfindahlÐHirschman_index}
%See the Brown and Warren-Boulton (1988) and Warren-Boulton (1990) texts cited below.
%}

%%%%%%%%%%%%%%%%%%%%%%%%
\subsection{ESS in political science}\label{ESSpolSect}
%%%%%%%%%%%%%%%%%%%%%%%%

In political science, ESS formulas have been used to set the effective number of parties in a political system. More precisely, the authors in \cite{Laakso79} proposed the effective number of parties using the following formula $\frac{1}{\sum_{n=1}^N {\bar w}_n^2}$ (the Hill number of order $2$ in ecology),
where $N$ is the total number of parties and ${\bar w}_n$ is the  proportion of votes of the $n$-th party.
%\url{https://en.wikipedia.org/wiki/Effective_number_of_parties}
%This measure is equivalent to the Herfindahl-Hirschman index, used in economics; the Simpson diversity index in ecology; the inverse participation ratio (IPR) in physics; and the RŽnyi entropy of order $\alpha=2$ in information theory. ${ }^{[8]}$
An alternative formula  was introduced in political science by Grigorii Golosov \cite{Golosov10}, 
%$$
%N=\sum_{i=1}^n \frac{p_i}{p_i+p_1^2-p_i^2}
\begin{align}\label{ESSGOL}
\mbox{ESS-Gol}_N({\bf {\bar w}})=\sum_{n=1}^N \frac{1}{1+\left(\frac{(\max \bar{w}_n)^2}{ \bar{w}_n}\right)- \bar{w}_n}=\sum_{n=1}^N \frac{\bar{w}_n}{\bar{w}_n+\left(\max \bar{w}_n \right)^2- \bar{w}_n^2},
\end{align}
that is also {\it proper and stable}. The value $\max \bar{w}_n$ denotes the portion of votes of the party that are obtained the greatest number of votes.   Other alternatives can be found in the literature \cite{Molinar91}.

%Here, $n$ is the number of parties, $p_i^2$ the square of each party's proportion of all votes or seats, and $p_1^2$ is the square of the largest party's proportion of all votes or seats.

%%%%%%%%%%%%%%%%%%%%%%
\subsection{ESS in quantum physics} 
\begin{figure}[!hbt]
\centering 
  \centerline{
\includegraphics[width=16cm]{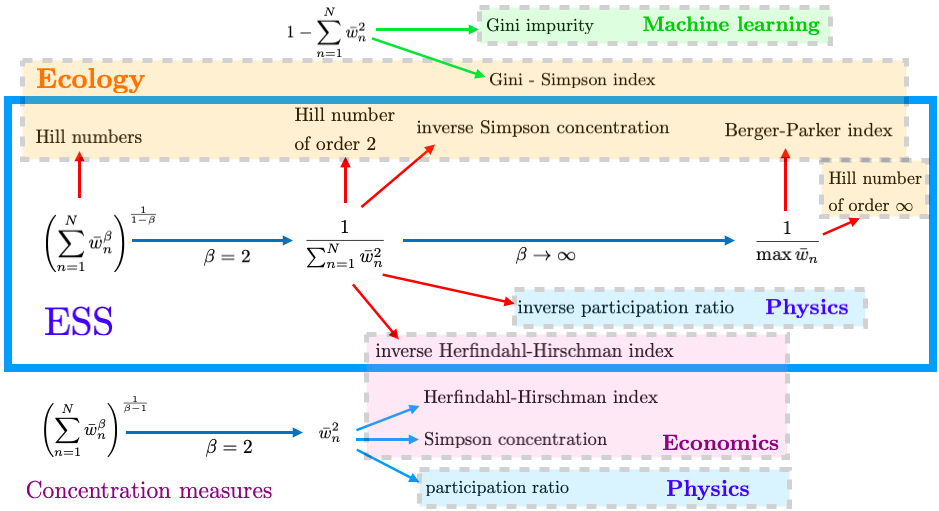}
   }
  \caption{{\small Graphical summary of the main nomenclature in different fields.   } }
  \label{allNames}
\end{figure}

%%%%%%%%%%%%%%%%%%%%%%
In quantum physics, there exists a quantity that is related to the formula
$\mbox{ESS-H}_N^{(2)}({\bf \bar{w}})=\sum_{n=1}^N {\bar w}_n^2$, 
that is called  {\it participation ratio} (PR) and the corresponding concentration $\mbox{Conc}_N^{(2)}({\bf \bar{w}})$ is known as
{\it  inverse participation ratio (IPR)}.   For a fully delocalized or spread state, we have the lowest value of the IPR, i.e., $\min (\mathrm{IPR})=1 / N$. On the other hand, for a fully localized state, we have the highest value of IPR, i.e., $\max (\mathrm{IPR})=1$.  IPR is also close to the concept of {\it purity} whereas PR is close to the concept of {\it separability}, employed both in quantum mechanics \cite{Zyczkowski98}. 
Moreover, since the purity $\mathcal{P}$ of a quantum state is a quantity such that $ \mathcal{P} \leq 1$, another concepts  naturally arises that is  {\it state mixedness} as the complement of purity, $\mathcal{M} = 1-\mathcal{P}$.  The quantum state is pure if $\mathcal{P}=1$. Figure \ref{allNames} summarizes the main nomenclature described so far.

%The smallest value of the IPR corresponds to a fully delocalized state, $\psi(x)=1 / \sqrt{N}$ for a system of size $N$, where the IPR yields $\sum_x|\psi(x)|^4=N /\left(N^{1 / 2}\right)^4=1 / N$. Values of the IPR close to 1 correspond to localized states (pure states in the analogy), as can be seen with the perfectly localized state $\psi(x)=\delta_{x, x_0}$, where the IPR yields $\sum_x|\psi(x)|^4=1$.

 %turns out to be useful. It is defined as the integral (or sum for finite system size) over the square of the density in some space, e.g., real space, momentum space, or even phase space, where the densities would be the square of the real space wave function $|\psi(x)|^2$, the square of the momentum space wave function $|\tilde{\psi}(k)|^2$, or some phase space density like the Husimi distribution, respectively. ${ }^{[6]}$

%%%%%%%%%%%%%%%%%%%%%%%%%%%%%%%%%%%%%%%%%%%%
\subsection{Application to model selection as effective number of components} \label{ESS-ENVsect}
%%%%%%%%%%%%%%%%%%%%%%%%%%%%%%%%%%%%%%%%%%%%

Model selection is a fundamental task in statistics and machine learning. An interesting scenario is when we have a family of nested models, where the model complexity can change since the number of parameters can vary (i.e., the dimension of the vector of parameters grows, building more complex models). The dimension of the vector of parameters is itself object of inference. This is the case of the order selection in polynomial regression problems or autoregressive schemes, variable selection, clustering, dimension reduction, just to name a few. 
In the literature, generally  cross-validation (CV) techniques \cite{bishop2007}, and information criteria \cite{Hannan79,SICpaper,schwarz1978estimating,Spiegelhalter02}  the procedures used to handle this problem.
More recently, other approaches  based on geometric considerations, have been also proposed in the literature, such as the automatic detection of an ``elbow'' or ``knee-point'' in a non-increasing curve describing a metric of performance of the model versus its complexity \cite{AEDpaperNuestro,Elbow2,Elbow3,Elbow4}. 
An effective number of variables/features (ENV) has been also proposed \cite{ENVpaper}. %Note that, {throughout} all the work, we use the words  {\it variables}, {\it components}, {\it features},   and/or {\it parameters} of a model as synonymous \cite{Elbow_paper,VarSelRev_paper}. 
The ENV index is inspired by the concept of maximum area-under-the-curve (AUC) in receiver operating characteristic (ROC) curves  \cite{HanleyAUC82} and the Gini inequality index,  described in Section \ref{other_ESS_Stable} and mentioned in Section \ref{ConcMeSect} %\cite{Lorenz05,Gini12,Gini13,Gini21}
\cite{Gini21}. In the variable selection scenario, the ENV index is given by 
 \begin{align}
I_{\mathrm{ENV}}=1+\frac{2}{V(0)} \sum_{k=1}^{N-1} V(k), \qquad \mbox{for  } V(0)\neq 0, \mbox{ and  }  V(N)=0,
 \end{align}
where $V(k)$ is a non-increasing error curve, e.g., the mean square errors (MSE), for a model that uses only $k\leq N$ input variables (instead of all the $N$ possible variables). By construction, it is always possible to have $V(N)=0$ (by a simple translation). It is possible to show that $1\leq I_{\mathrm{ENV}} \leq N$.  The ENV index could be defined also for non-decreasing curve by the alternative definition
  \begin{align}\label{IncV}
I_{\mathrm{ENV}}=1+\frac{2}{V(N)} \sum_{k=1}^{N-1} V(k), \qquad \mbox{for  } V(N)\neq 0, \mbox{ and  }  V(0)=0.
 \end{align}
Thus, we can convert the  ENV index in an ESS formula building the curve $V(k)$ as follows:
 \begin{itemize}
\item Sort in ascending order the normalized weights as
$$
\bar{w}_{(1)} \leq \bar{w}_{(2)} \leq \ldots \leq \bar{w}_{(N)}.
$$
\item Build a non-decreasing curve $V(k)$, as in \eqref{IncV}, following the recursion:
\begin{align}
V(k)=\sum_{i=1}^k \bar{w}_{(i)}=V(k-1)+\bar{w}_{(k)},
\end{align}
starting with $V(0)=0$. Note that we always have $V(N)=1$.
 \end{itemize}
The corresponding ESS formula is 
\begin{align}
\mbox{ESS-ENV}_N({\bf \bar{w}})=1+2 \sum_{k=1}^{N-1} \sum_{i=1}^k \bar{w}_{(i)}.
\end{align}
Note that  $\mbox{ESS-ENV}_N({\bf \bar{w}}^{(j)})=1+0=1$ for all $j$, and
\begin{align}
\mbox{ESS-ENV}_N({\bf \bar{w}}^*)=1+\frac{2}{N} \sum_{k=1}^{N-1} k=1+\frac{2}{N} \frac{(N-1) (N)}{2}=1+N-1=N.
\end{align}

{\rem It is possible to show that $\mbox{ESS-ENV}_N({\bf \bar{w}})$ is proper and stable. Furthermore, it 
coincides with $\mbox{ESS-Gini}({\bf \bar{w}})$ in Eq. \eqref{GINIess}, i.e.,  $\mbox{ESS-ENV}_N({\bf \bar{w}})=\mbox{ESS-Gini}_N({\bf \bar{w}})$. Recall that there exist different formulations of the Gini coefficient \cite{Gini13}. The closest one in this framework is related to the Lorenz curve \cite{Lorenz05}.}

{\rem This section opens the possibility to apply the ESS formulas as effective number fo components in model selection problems. Indeed, given a  non-increasing  error curve $V(k)$, i.e., $V(k-1)\leq V(k)$, we can build the normalized weights in this way:
\begin{align}\label{Appl_to_ENC}
d_k=V(k-1)-V(k), \qquad \bar{w}_k=\frac{d_k}{\sum_{i=1}^N d_i},  
\end{align}
for all $k=1,...,N$. Then the ESS formula can be applied to the vector ${\bf \bar{w}}=[\bar{w}_1,...,\bar{w}_N]$.}

%$$
%I_{\mathrm{ENV}}=1+2 \sum_{k=2}^{N-1} \frac{\bar{w}_{(k+1)}-\bar{w}_{(k)}}{\bar{w}_{(2)}-\bar{w}_{(1)}},  \quad \bar{w}_{(2)}-\bar{w}_{(1)} \neq 0
%$$

%%%%%%%%%%%%%%%%%%
\section{Numerical experiments}
%%%%%%%%%%%%%%%%%%

\subsection{Analyzing the Huggins-Roy family }

Since all the ESS functions in the Huggins-Roy family are proper and stable and, since this family contains the main relevant formulas, we focus the numerical experiments on this family. First of all, we recall the theoretical definition of ESS in Eq. \eqref{DEF_ESS_1},
\begin{equation}
\mbox{ESS}_{\texttt{teo}}(h)=N\frac{\mbox{var}_\pi[{\widehat I}]}{\mbox{var}_q[{\widetilde I}]}.
\end{equation}
where, for simplicity, we consider a scalar $x\in \mathbb{R}$ the use of the integrand  $h( x)=x$ (in the definition above, we have clarified the dependence on the function $h$). Namely, ${\widehat I}$ and ${\widetilde I}$ are estimators of the expected value of a random variable $X$ with a target pdf $\bar{\pi}(x)$ (defined below). In this numerical example, we compute approximately via Monte Carlo the theoretical definition $\mbox{ESS}_{\texttt{teo}}$, and compare them with the G-ESS functions $\mbox{ESS-H}_N^{(\beta)}$. More specifically, we consider a univariate standard Gaussian density as target pdf, 
\begin{equation}
\bar{\pi}(x)=\mathcal{N}(x;0,1),
\end{equation}
and also a Gaussian proposal pdf,
\begin{equation}
q(x)=\mathcal{N}(x;\mu_p,\sigma_p^2),
\end{equation}
with mean $\mu_p$ and variance $\sigma_p^2$. In all the experiments, we consider $N=1000$.

%%%%%%%%%%%%%%%%%%%%%%%%%%%%%%%%%
\subsubsection{Varying the proposal mean $\mu_p$}
\begin{figure}[!hbt]
\centering 
  \centerline{
  \subfigure[\label{fig0simu_a}]{\includegraphics[width=6cm]{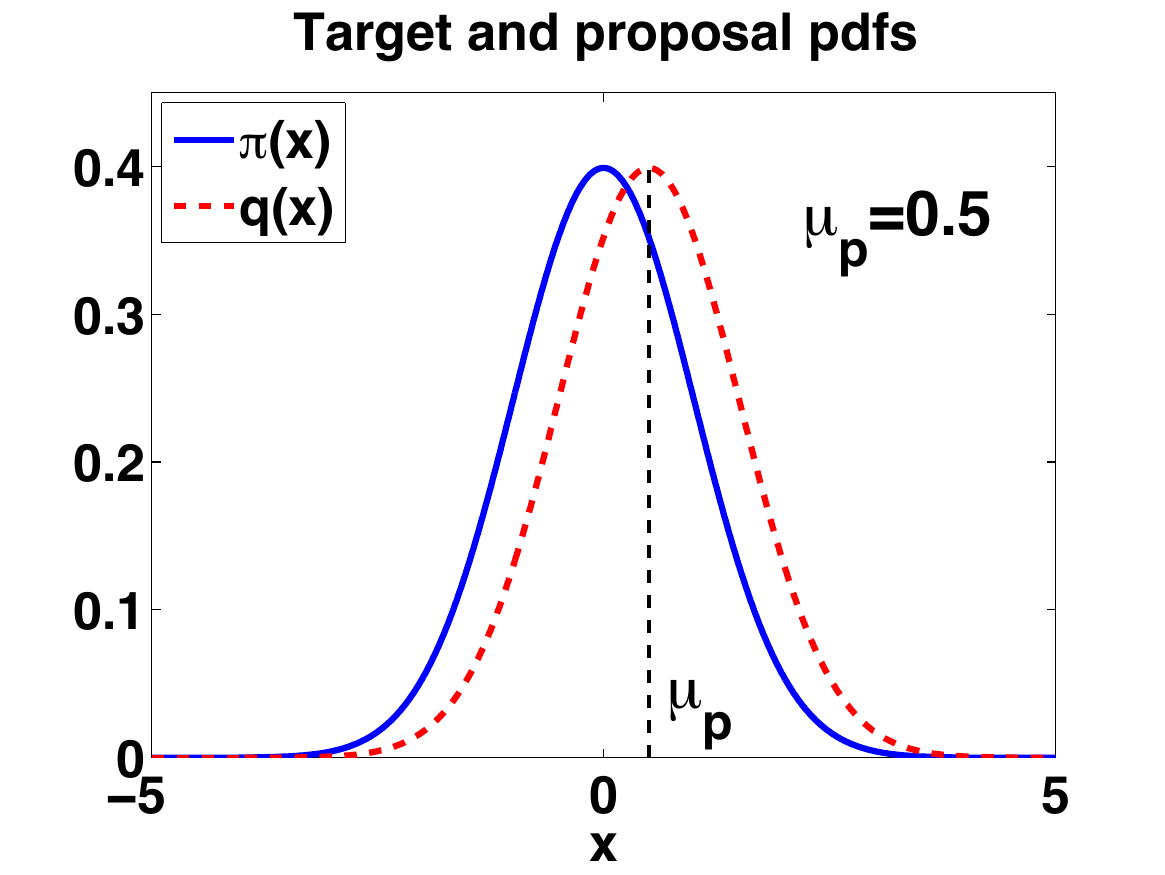}}
  \subfigure[\label{fig0simu_b}]{\includegraphics[width=6cm]{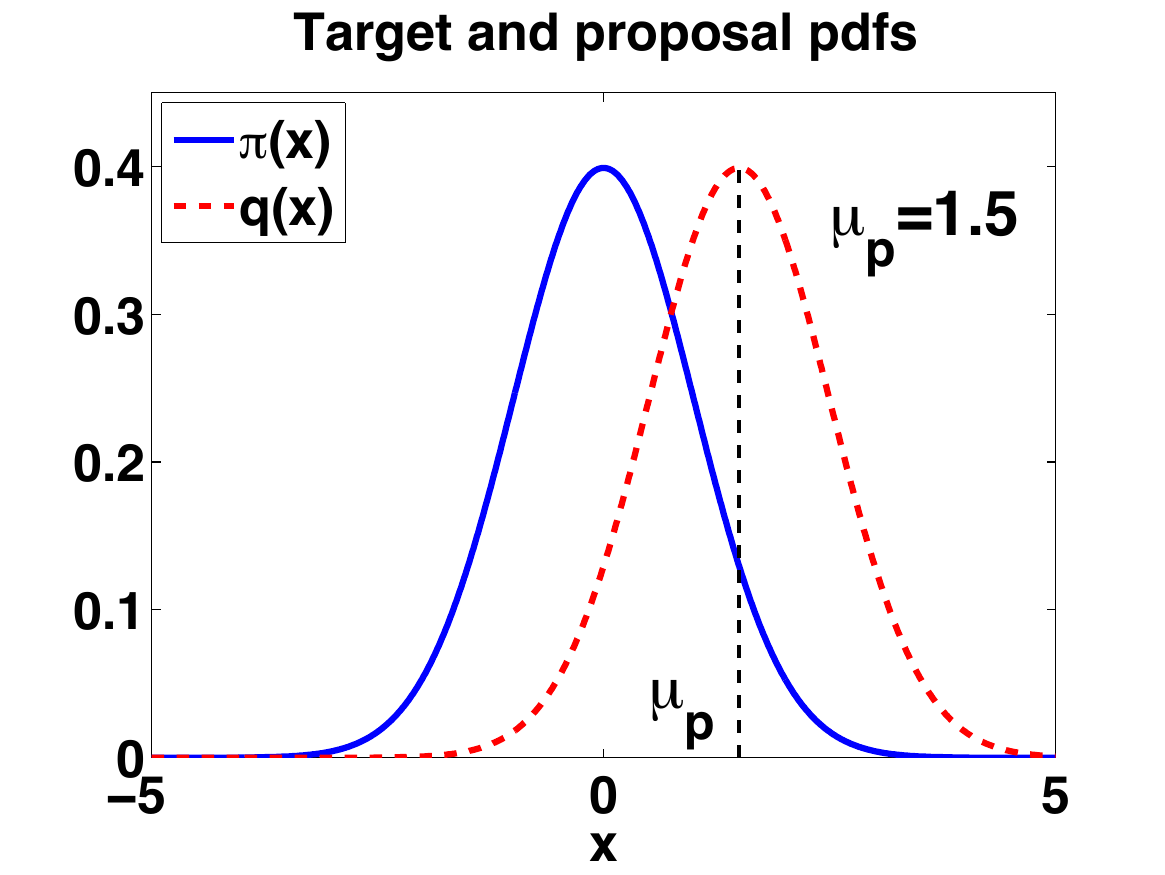}}
   \subfigure[\label{fig0simu_c}]{\includegraphics[width=6.1cm]{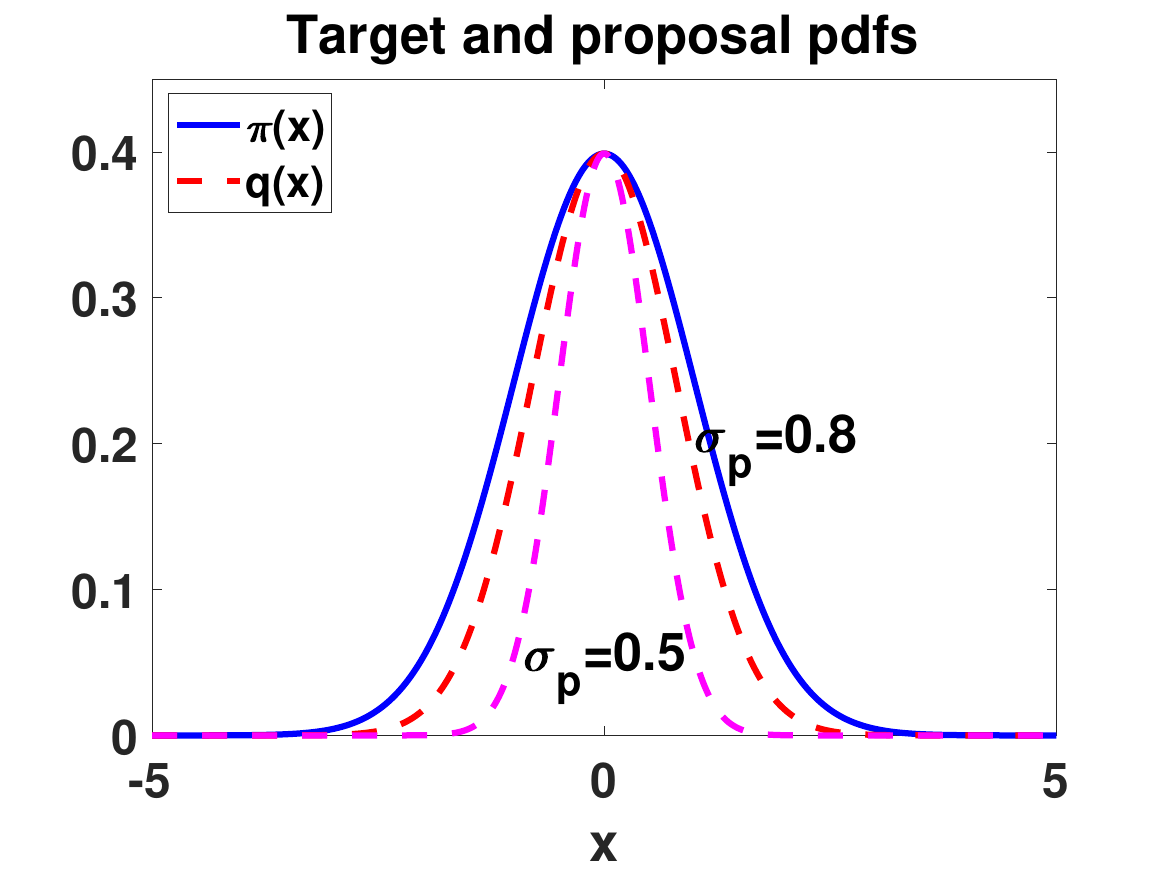}}
 %%%%  \subfigure[]{\includegraphics[width=6cm]{FigSIMU_varVar_v2.pdf}}
   }
  \caption{{\small Target and proposal pdfs: {\bf (a)}-{\bf (b)} with $\mu_p\in\{0.5,1.5\}$. The variances in both is set to $1$. {\bf (c)} here $\mu_p=0$ and $\sigma_p\in\{0.5,0.8\}$. } }
\label{fig0simu}
\end{figure}

\begin{figure}[!hbt]
\centering 
  \centerline{
  \subfigure[\label{fig1simu}]{\includegraphics[width=9.5cm]{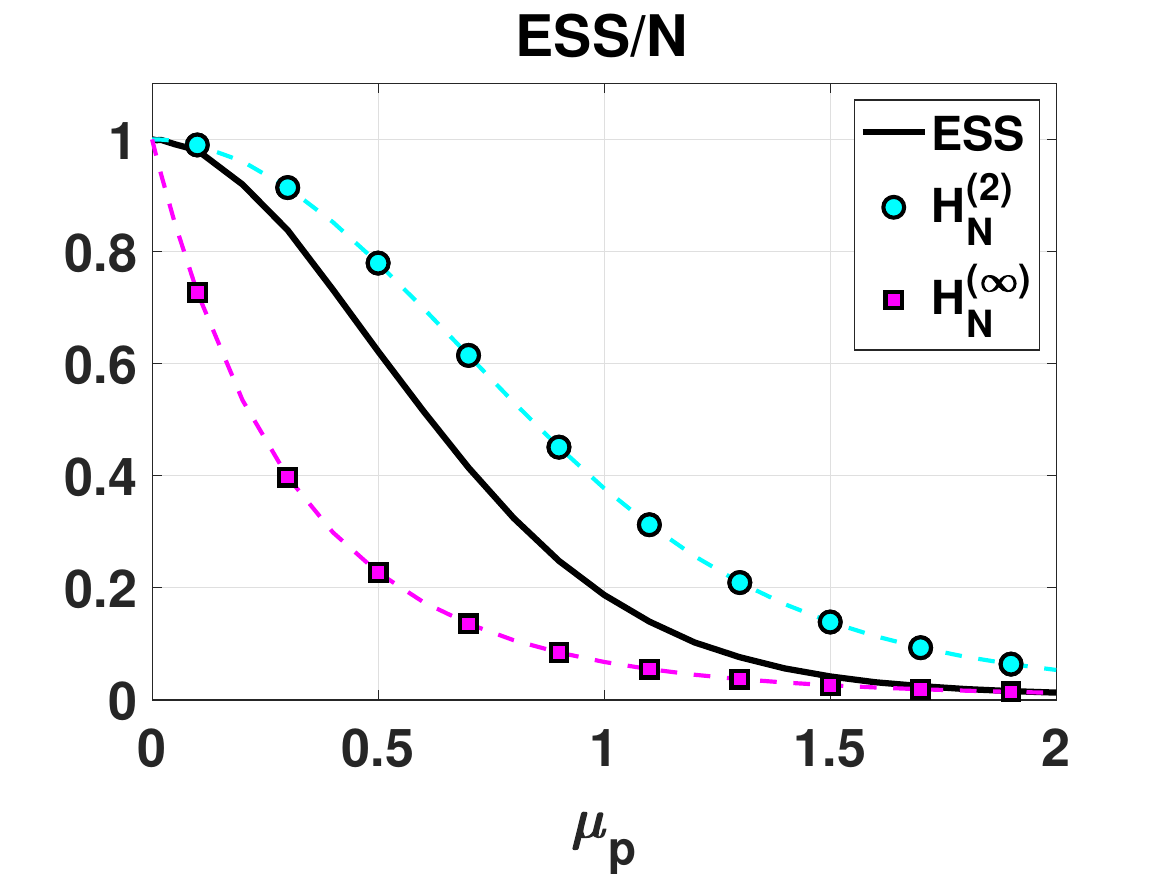}}
  \subfigure[\label{fig2simu}]{  \includegraphics[width=9.5cm]{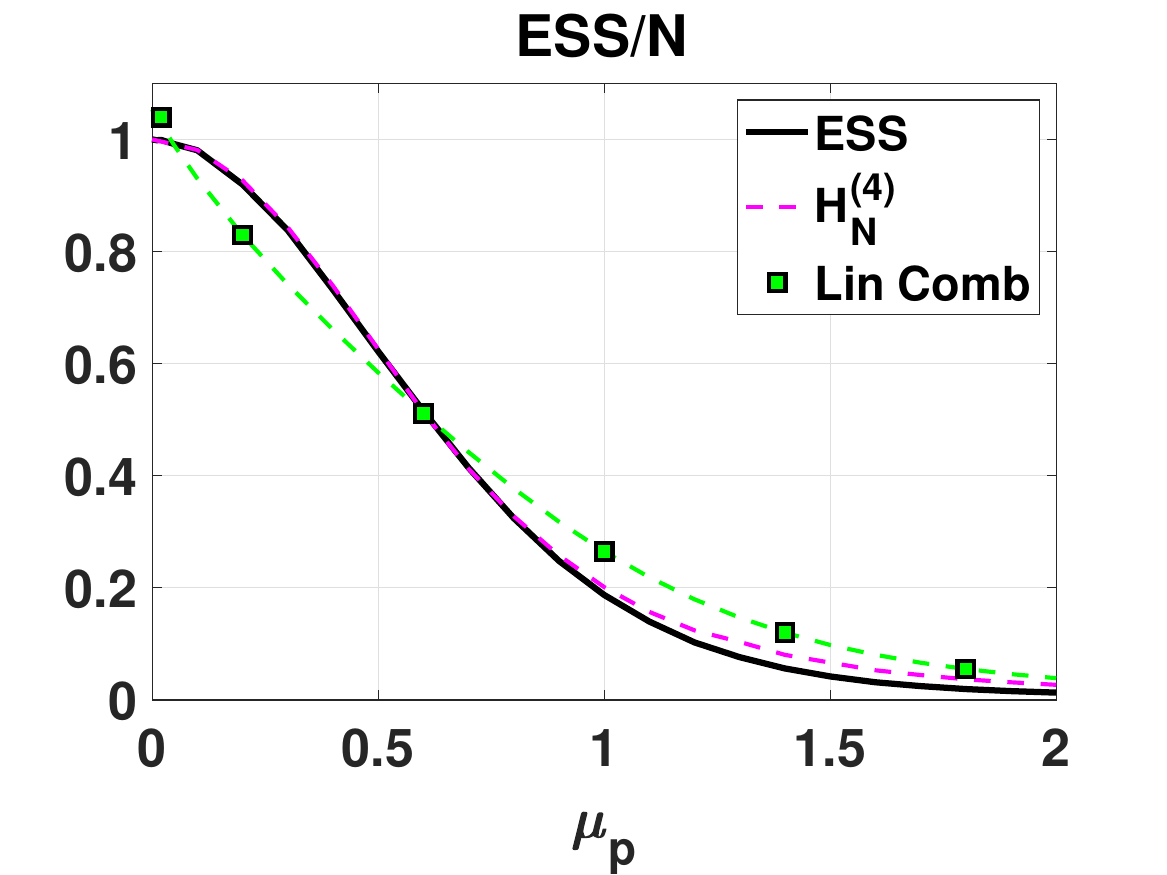}}
   }
  \caption{{\small Ratio of ESS values over $N$ (with $N=1000$) versus $\mu_p$. The curve corresponding to  theoretical ESS value, i.e.,  $\mbox{ESS}_{\texttt{teo}}/N$ is shown in black solid line in both figures. In {\bf (a)} the curves of $\mbox{ESS-H}_N^{(2)}/N$ (circles) and $\mbox{ESS-H}_N^{(\infty)}/N$ (squares) are also depicted. In {\bf (b)} we show the curves $\mbox{ESS-H}_N^{(4)}/N$ (dashed line) and the linear combination in Eq. \eqref{LinC}-\eqref{LinC2} (squares), as well. The approximation provided by $\mbox{ESS-H}_N^{(4)}$ is virtually perfect for $\mu_p\leq 1$.  } }
\end{figure}

%%%%%%%%%%%%%%%%%%%%%%%%%%%%%%%%%
  In a first analysis, we keep fixed $\sigma_p=1$ and vary $\mu_p\in[0,2]$. Figures \ref{fig0simu_a}-\ref{fig0simu_b} depict two scenarios in this experimental setup, corresponding to two specific values of $\mu_p$, $0.5$ and $1.5$.  Clearly, for $\mu_p=0$ we have the ideal Monte Carlo case, $q(x)\equiv \bar{\pi}(x)$. As $\mu_p$ increases, the proposal becomes more different from $\bar {\pi}$.  We recall that $N=1000$. %Figure \ref{fig1simu} shows the (approximated) theoretical ESS curves and the curves corresponding to  different ESS formulas, averaged over $10^5$ independent runs. More specifically, we provide the rates $\frac{ESS}{N}$. Note that $\frac{1}{N}\leq\frac{ESS}{N}\leq 1$.
%whereas  for $N=5$ and S2 we only provide $ESS_{MSE}$ since the bias is big for small value of $\sigma_p$ so that it is difficult to obtain reasonable and meaningful values of $ESS_{var}$.
Figure \ref{fig1simu} shows the  theoretical $\mbox{ESS}_{\texttt{teo}}/N$ curves (solid line)  $\mbox{ESS-H}_N^{(2)}/N$ (circles) and $\mbox{ESS-H}_N^{(\infty)}/N$ (squares),  averaged over $10^5$ independent runs. Note that $\frac{1}{N}\leq\frac{ESS}{N}\leq 1$. 

 \paragraph{Optimal linear combination of $\mbox{ESS-H}_N^{(2)}$ and $\mbox{ESS-H}_N^{(\infty)}$.}  The functions $\mbox{ESS-H}_N^{(2)}$ and $\mbox{ESS-H}_N^{(\infty)}$ are the most 
 used and suggested formulas in different studies \cite{Huggins15,ESSarxiv16}. Moreover, at least in this simulation scenario, they seem to play the role of upper bound and lower bound of the true value, as shown by Figure \ref{fig1simu}. For this reason, we also consider the linear combination of the G-ESS formulas $\mbox{ESS-H}_N^{(2)}$ and $\mbox{ESS-H}_N^{(\infty)}$,
\begin{equation}
\label{LinC}
\mbox{Comb-ESS}_N({\bf \bar w})= a_1  \mbox{ESS-H}_N^{(2)}({\bf \bar w})+  a_2  \mbox{ESS-H}_N^{(\infty)}({\bf \bar w}).
\end{equation}
This example suggests the use of 
\begin{eqnarray}
\label{LinC2_coeff}
&&a_1=0.6245, \nonumber \\
&& a_2=0.4289,
\end{eqnarray}
 obtained using a Least Squares (LS) regression in order to obtain an expression $\mbox{Comb-ESS}_N({\bf \bar w})$ as close as possible to the theoretical ESS curve. 
 
 \paragraph{Optimal $\beta$ for $\mbox{ESS-H}_N^{(\beta)}({\bf \bar w})$.} Furthermore, we have computed the curves (as function $\beta$) of $\mbox{ESS-H}_N^{(\beta)}({\bf \bar w})$ for different values of $\beta$, considering a thin grid of $\beta$ values from $0.2$ to $50$ with a step of $0.01$ (i.e., $\beta \in \mathcal{G}$ denoting $\mathcal{G}$ the thin grid).
  We consider a $L_1$ distance between each $\mbox{ESS-H}_N^{(\beta)}({\bf \bar w})$ curve and the theoretical ESS curve,\footnote{Recall that these curves are functions of $\mu_p$ and are averaged over $10^5$ independent runs.}, i.e., $|\mbox{ESS-H}_N^{(\beta)}-\mbox{ESS}_{\texttt{teo}}|$,  and compute
\begin{equation}
\beta^*=\arg \min_{\beta \in \mathcal{G}} |\mbox{ESS-H}_N^{(\beta)}-\mbox{ESS}_{\texttt{teo}}|.
\end{equation}
With this procedure, we obtain 
$$
\beta^*\approx 4.
$$
 \paragraph{Discussion of the results.} Figure \ref{fig2simu} shows the curves of the ESS rates corresponding to the theoretical ESS curve (solid line), the best linear combination corresponding to the Eqs. \eqref{LinC}-\eqref{LinC2} (squares) and the curve corresponding to $\mbox{ESS-H}_N^{(\beta^*)}$ (dashed line). First of all, we can note that the linear combination can return values greater than 1 (recall that we are considering $\mbox{ESS}/N$).
 Moreover, we can see that the curve corresponding to $\mbox{ESS-H}_N^{(4)}({\bf \bar{w}})$ fits particularly well in this numerical setup, providing a very close to the theoretical ESS curve. Observe that the approximation provided by $\mbox{ESS-H}_N^{(4)}$ is virtually perfect for $\mu_p\leq 1$. Hence, in this kind of scenario, we would suggest the use of the expression
 \begin{eqnarray}
\mbox{ESS-H}_N^{(4)}({\bf \bar{w}})&=&\left(\frac{1}{\sum_{n=1}^N {\bar w}_n^4}\right)^{\frac{1}{3}}. Ê
\end{eqnarray}

\begin{figure}[!hbt]
\centering 
  \centerline{
  \subfigure[\label{fig3simu}]{\includegraphics[width=9.5cm]{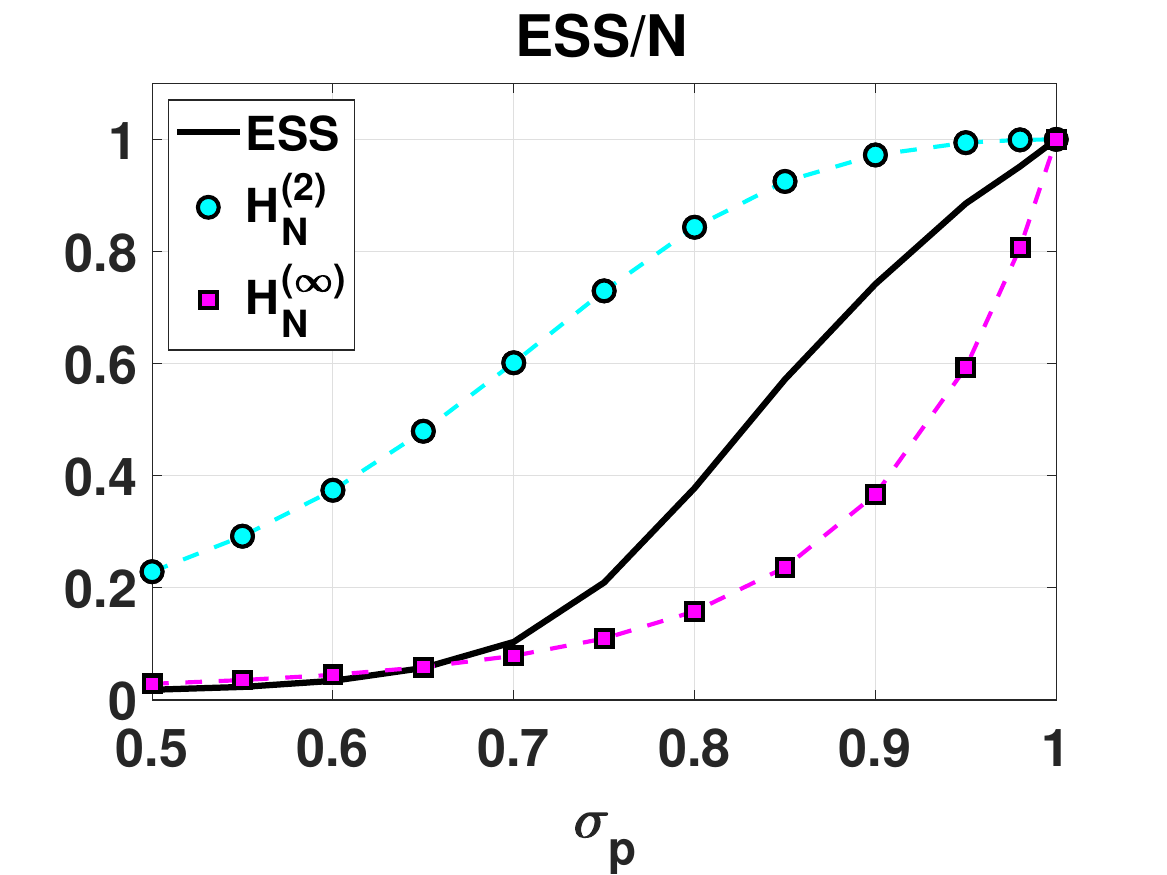}}
  \subfigure[\label{fig4simu}]{  \includegraphics[width=9.5cm]{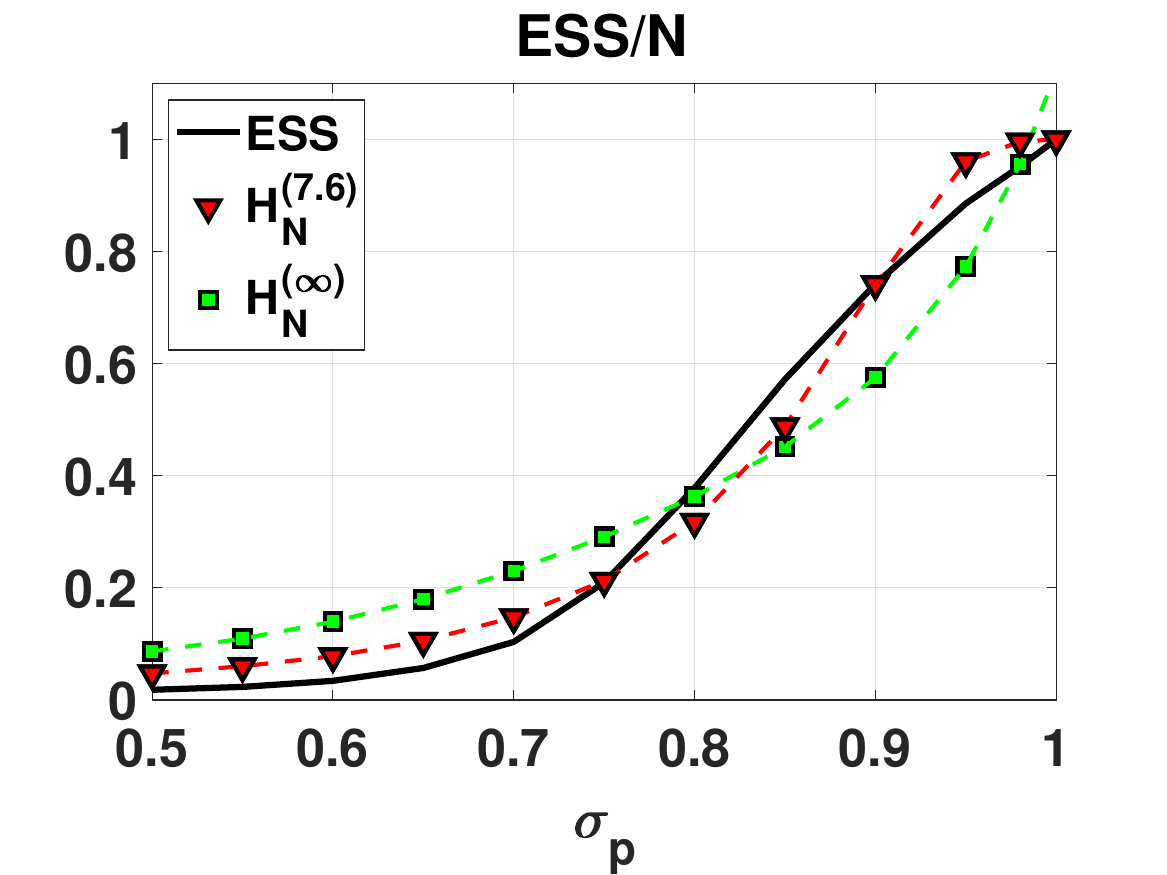}}
   }
  \caption{{\small Ratio of ESS values over $N$ (with $N=1000$) versus $\sigma_p$. The curve corresponding to  theoretical ESS value, i.e.,  $\mbox{ESS}_{\texttt{teo}}/N$ is shown in black solid line in both figures. In {\bf (a)} the curves of $\mbox{ESS-H}_N^{(2)}/N$ (circles) and $\mbox{ESS-H}_N^{(\infty)}/N$ (squares) are also depicted. In {\bf (b)} we show the curves $\mbox{ESS-H}_N^{(7.6)}/N$ (dashed line) and the linear combination in Eq. \eqref{LinC3} (squares), as well.   } }
\end{figure}

%%%%%%%%%%%%%%%%%%%%%%%%%%%%%%%%%
\subsubsection{Varying the proposal standard deviation $\sigma_p$}
%%%%%%%%%%%%%%%%%%%%%%%%%%%%%%%%%
Now, we keep fixed $\mu_p=0$ and vary the standard deviation of the proposal $\sigma_p\in[0.5,1]$. Figure \ref{fig0simu_c} depicts the target density and the proposal density for two specific values of $\sigma_p$, $0.5$ and $0.8$, used in  this experimental setup.   
We recall that $N=1000$ and the results have been averaged over $10^5$ independent runs. In Figure \ref{fig3simu}, we can observe the results of  $\mbox{ESS}_{\texttt{teo}}/N$ versus $\sigma_p$ (in solid line),  jointly with the curves $\mbox{ESS-H}_N^{(2)}/N$ (given with circles) and $\mbox{ESS-H}_N^{(\infty)}/N$ (shown with squares).

 \paragraph{Optimal linear combination of $\mbox{ESS-H}_N^{(2)}$ and $\mbox{ESS-H}_N^{(\infty)}$.}  Since the formulas $\mbox{ESS-H}_N^{(2)}$ and $\mbox{ESS-H}_N^{(\infty)}$ are the most 
 used in practice, again we consider the linear combination of the G-ESS formulas $\mbox{ESS-H}_N^{(2)}$ and $\mbox{ESS-H}_N^{(\infty)}$,
\begin{equation}
\label{LinC2}
\mbox{Comb-ESS}_N({\bf \bar w})= a_1  \mbox{ESS-H}_N^{(2)}({\bf \bar w})+  a_2  \mbox{ESS-H}_N^{(\infty)}({\bf \bar w}),
\end{equation}
where in this scenario we get by LS solution
\begin{eqnarray}
\label{LinC3}
&&a_1=0.2715, \nonumber \\
&& a_2=0.8483,
\end{eqnarray}
hence $\mbox{ESS-H}_N^{(\infty)}$ takes more importance in this scenario.
Figure \ref{fig4simu} provides the curve corresponding to $\mbox{Comb-ESS}_N({\bf \bar w})/N$ with a dashed line and green squares.
 
 \paragraph{Optimal $\beta$ for $\mbox{ESS-H}_N^{(\beta)}({\bf \bar w})$.} Furthermore, we have computed the curves (as function $\beta$) of $\mbox{ESS-H}_N^{(\beta)}({\bf \bar w})$ for different values of $\beta$, considering a  grid of values of $\beta$ denoted as $\mathcal{G}$.
  We consider a $L_1$ distance between each $\mbox{ESS-H}_N^{(\beta)}({\bf \bar w})$ curve and the theoretical ESS curve, and compute
\begin{equation}
\beta^*=\arg \min_{\beta \in \mathcal{G}} |\mbox{ESS-H}_N^{(\beta)}-\mbox{ESS}_{\texttt{teo}}|.
\end{equation}
In this scenario, we obtain 
$$
\beta^*\approx 7.6.
$$
The corresponding curve is depicted in Figure \ref{fig4simu} with a dashed line and red  triangles. We can see that we obtain a very good approximation of $\mbox{ESS}_{\texttt{teo}}/N$, but slightly worse than in the case described in the previous section. Moreover, here the optimal $\beta^*$ is $\approx 7.6$ whereas, in the previous section, was $\beta^*$ is $\approx 4$. 

 \paragraph{Discussion of the results.} Figure \ref{fig4simu} shows the curves of the ESS rates corresponding to the theoretical ESS curve (solid line), the best linear combination corresponding to the Eqs. \eqref{LinC}-\eqref{LinC2} (green squares) and the curve corresponding to $\mbox{ESS-H}_N^{(\beta^*)}$ (red triangles). Again the linear combination can return values greater than 1 (recall that we are considering $\mbox{ESS}/N$). This behavior  could be exploited in future works since actually $\mbox{ESS}_{\texttt{teo}}/N$ can exceed $1$ (see  \cite[Section 3.3]{elvira2018rethinking}).
 Moreover, we can see that $\mbox{ESS-H}_N^{(7.6)}({\bf \bar{w}})$ performs particularly well in this scenario, providing a close to the theoretical ESS curve. Hence, in this setup, we would suggest the use of  $\mbox{ESS-H}_N^{(7.6)}({\bf \bar{w}})$. Only for simplicity in computation and comparison, one could consider the closest integer and use $\beta=8$, 
  \begin{eqnarray}
\mbox{ESS-H}_N^{(8)}({\bf \bar{w}})&=&\left(\frac{1}{\sum_{n=1}^N {\bar w}_n^8}\right)^{\frac{1}{7}}. Ê
\end{eqnarray}
 Finally, it is important to remark that even if the optimal $\beta^*\approx 7.6$ (or $8$)  is different from the value $\beta^*\approx 4$ suggested in the previous section, however both values differ from $2$ (that corresponds to the typical formula employed in the literature) and  both values are bigger than $2$. The expression with $\beta \rightarrow \infty$, i.e., $\mbox{ESS-H}_N^{(\infty)}=\frac{1}{\max \bar{w}_n}$ seems that can be employed as a lower bound for the theoretical value $\mbox{ESS}_{\texttt{teo}}$, in both setups.
  These considerations can be  relevant clues for future applications and studies.

 %%%%%%%%%%%%%%%%%%%%%%%%%%%%%%%%%%%%%
\subsection{Application to a variable selection in a regression problem with real data}\label{Ex2Sound}
%%%%%%%%%%%%%%%%%%%%%%%%%%%%%%%%%%%%%

Finding the connections with other fields creates the opportunities for new applications for the ESS formulas. As we described in Section \ref{ESS-ENVsect}, the ESS can be applied in a feature selection problem to find the effective number of components. In this section, we provide an example of this application with a real dataset.
\newline
Let us consider a regression problem, where we observe a dataset of $N$ pairs $\left\{\mathbf{x}_n, y_n\right\}_{n=1}^N$, with each input vector $\mathbf{x}_n=\left[x_{n, 1}, \ldots, x_{n, K}\right]$ is formed by $K$ variables, and the outputs $y_n$ 's are scalar values \cite{OurPaperSound}. We consider the case that being $K \leq N$ and assume a linear observation model,
\begin{align}\label{ModelExpSoundScape}
y_n=\theta_0+\theta_1 x_{n, 1}+\theta_2 x_{n, 2}+\ldots \theta_K x_{n, K}+\epsilon_n,
\end{align}
where $\epsilon_n$ is a Gaussian noise with zero mean and variance $\sigma_\epsilon^2$, i.e., $\epsilon_n \sim \mathcal{N}\left(\epsilon | 0,\sigma_\epsilon^2 \right)$. More specifically, in this real dataset \cite{OurPaperSound,san2024variable,fan2017emo}, we have $K=122$ features and $N=1214$ number of data points ${\bf x}_i$.  
We focus on the first of the two outputs in the dataset (called ``arousal'').
We set $V(k)=-2 \log \left(\ell_{\max }\right)$ with $\ell_{\max }=\max _\theta p\left(\mathbf{y} | \boldsymbol{\theta}_k\right)$ with $k \leq K$, after ranking the 122 variables (see \cite{OurPaperSound}), where the likelihood function $p\left(\mathbf{y} | \boldsymbol{\theta}_k\right)$ is induced by the Eq. \eqref{ModelExpSoundScape}. In order to find the effective number of variables $N_{\texttt{eff}}\leq K=122$, we compare with different well-known information criteria\footnote{Considering the cost function $C(k)=V(k)+\lambda k$, each  information criterion suggests the use of a different parameter $\lambda$.},  AIC, BIC and HQIC, and other methods provided in the literature.
For the spectral information criterion (SIC), we test two confidence internal parameter to $95\%$ and $99\%$. We also test different stable ESS formulas obtained the weights as in Eq. \eqref{Appl_to_ENC}. 
We test the expressions in the Huggins-Roy family, $\mbox{ESS-H}_N^{(\beta)}$, with $\beta \rightarrow 1$, $\beta=2$, $\beta \rightarrow \infty$, and the other stable formulas given in the Eqs. \eqref{Nplus},  \eqref{Qplus},  \eqref{GINIess}, and  \eqref{ESSGOL}.
All the results are {\it rounded} to the closest integer. Thus, the results provided by each method are given in Table 
\ref{Table1}.

\begin{table}[!h]
\caption{Results in the variable selection example with a real dataset. } \label{Table1}
\begin{center}
%\footnotesize
%\scriptsize
\begin{tabular}{|c|c|c|c|c|c|c|c|}
\hline
{\bf Scheme} & AIC & BIC & HQIC & UAED & SIC-95 & SIC-99  & ENV \\
\hline
$N_{\texttt{eff}}$  & 44 & 17 &  41 &  11 & 7 & 17  &  13  \\
 %&  &  &   &   &  &   &  {\bf 12-13}  \\
\hline
{\bf Ref.}  & \cite{Spiegelhalter02} &    \cite{schwarz1978estimating} &   \cite{Hannan79}    &   \cite{AEDpaperNuestro}    &  \cite{SICpaper}  &  \cite{SICpaper} &  \cite{ENVpaper} \\ 
\hline
\multicolumn{8}{c}{ } \\
\hline
{\bf ESS formula} &  $\beta \rightarrow 1$  & $\beta=2$Ê&  $\beta \rightarrow \infty$ &  Plus   &  Q  & Gini  & Gol \\
\hline
$N_{\texttt{eff}}$  & 10 &  5  &  3   &  11     &  24   & 11 &  4 \\
 %&  &  &   &   &  &   &  {\bf 12-13}  \\
\hline
{\bf Eq.}  &\eqref{PerpEq}&  \eqref{FirstDef_P} &  \eqref{MaxEq}    &  \eqref{Nplus}     &  \eqref{Qplus} &  \eqref{GINIess} &  \eqref{ESSGOL} \\ 
\hline
\end{tabular}
\end{center} 
\end{table}

\noindent
After an exhaustive analysis, the authors in \cite[Section 4-C]{OurPaperSound} suggest that there are 7 very relevant variables  (level 1 of  \cite[Section 4-C]{OurPaperSound}),  other 7 relevant variables (level 2 of  \cite[Section 4-C]{OurPaperSound}) and other 2 variables in a level 3 of importance  \cite[Section 4-C]{OurPaperSound}, hence, overall 16 variables among very relevant, relevant, and important ones (16 over 122 possible features).  The minimum value in Table \ref{Table1} is 3, provided by $\mbox{ESS-H}_N^{(\infty)}$, whereas the maximum value is 44 given by AIC. These values and the rest of results in Table \ref{Table1} are in line with the conclusions in \cite{OurPaperSound}.
More specifically,  the results given by the SIC-99, BIC,  UAED,  ENV, the perplexity $\mbox{ESS-H}_N^{(1)}$, $\mbox{ESS-Plus}$ and $\mbox{ESS-Gini}$ are $10 \leq N_{\texttt{eff}}\leq 17$, and are close to the results of the analysis in \cite{OurPaperSound}. Hence, in this experiment, some ESS formulas like the perplexity $\mbox{ESS-H}_N^{(1)}$, $\mbox{ESS-Plus}$ and $\mbox{ESS-Gini}$, seem to provide good performance as effective number of components in model selection.

\section{Conclusions}
%\vspace{-0.3cm}
 %%%%%%%%%%%%%%%%%% %%%%%%%%%%%%%%%%%%
 %%%%%%%%%%%%%%%%%% %%%%%%%%%%%%%%%%%%
 %%%%%%%%%%%%%%%%%% %%%%%%%%%%%%%%%%%%
 %%%%%%%%%%%%%%%%%% %%%%%%%%%%%%%%%%%%
 In this work, we have analyzed alternative effective sample size (ESS) measures for Monte Carlo algorithms based on the importance sampling techniques. We have remarked the connection to the practical ESS formulas used in the literature and  entropy families \cite{Cover91}.
 We have shown that all the ESS functions included in the Huggins-Roy's ESS family fulfill all the required theoretical conditions described in \cite{ESSarxiv16},  and we have also highlighted the relationship of this family with the R\'enyi entropy \cite{Cover91}. We have also shown the application of the Gini impurity index as ESS formula and its connection to the Tsallis entropy. 
 \newline
 Furthermore, we have studied the performance of different Huggins-Roy's ESS formulas  by numerical simulations, introducing also an optimal linear combination of the most promising ESS indices. In two numerical examples, we have obtained the best ESS approximations within the Huggins-Roy's family in two different setups, $\mbox{ESS}=\left(\frac{1}{\sum_{n=1}^M {\bar w}_n^{4}}\right)^{1/3}$ and  $\mbox{ESS}=\left(\frac{1}{\sum_{n=1}^M {\bar w}_n^{8}}\right)^{1/7}$. These formulas provide a good approximation (and in the first case almost a perfect match) of the theoretical ESS values, in two different considered experimental scenarios. Moreover, the expression  $\mbox{ESS}=\frac{1}{\max \bar{w}_n}$, which corresponds to  $\beta \rightarrow \infty$, also provides good performance in some specific cases (and playing the role of a lower bound of the ESS measures in other cases). All these considerations suggest us that the use of a  $\beta>2$ can more adequate in practical applications, e.g.,  in order to fight the sample degeneracy and impoverishment within a particle filtering algorithm.
 \newline
The relationship with the entropy families has also clarified the connections with other fields: possible applications in ecology, economics, political science, and machine learning have been discussed.  The application of the ESS expressions as the effective number of components in model selection seems to be promising but should be further investigated and tested. Moreover, the construction of these connections with other fields can also yield novel contributions in the IS context. As a final consideration, finding a novel and broader family that contains all the stable ESS formulas (that do not belong to the Huggins-Roy's family) could be object of future research. 
 
{\small
\section*{Acknowledgement}

 The work was partially supported  by the project Starting Grant for Rttb, BA-GRAPH ``Efficient Bayesian inference for graph-supported data", of the University of Catania (UPB-28722052144),  by the project LikeFree-BA-GRAPH funded by ``PIAno di inCEntivi per la RIcerca di Ateneo 2024/2026'' of the University of Catania (UPB-28722052159), Italy.}

%\section*{Acknowledgement}
%Luca Martino acknowledges support by the Agencia Estatal de Investigaci{\'o}n AEI (project
%SPGRAPH, ref. num. PID2019- 105032GB-I00).
%This work has been supported by the Grant 2014/23160-6 of S\~ao Paulo Research Foundation (FAPESP). %and by the Grant 305361/2013-3 of National Council for Scientific and Technological Development (CNPq).

%\bibliographystyle{elsarticle-num} 
%\bibliographystyle{IEEEtran}
%\bibliography{bibliografia,biblioFading}
%\bibliography{bib_battery,other,traffic_new,bibliografia}
\bibliographystyle{plain}
\bibliography{bibliografia}

%%%%%%%%
\appendix
%%%%%%%%
%%%%%%%%%%%%%%%%%%%%%%%%%%%%%%%%%%%%%%%
\section{Probabilistic interpretation}
\label{App1}
%%%%%%%%%%%%%%%%%%%%%%%%%%%%%%%%%%%%%%%
 Let us define a pair of random variables $\{X_t,Z_t\}$ that corresponds to generate  a random pair of samples that are independently drawn with replacement according to the pmf defined by ${\bar w}_n$, with $n=1,...,N$. We denote with $t\in \mathbb{N}$  a sub-index corresponding the trial/experiment. We perform different {\it independent} trials. Let also define the random variable
 $$
 T=\{\mbox{minimum $t\in  \mathbb{N}$ such that $X_t=Z_t$} \}.
 $$
 We aim to compute the expected number of trials needed to obtain a {\it first} pair containing the same sample twice, i.e., 
 \begin{align}\label{ETformula}
 E[T]=\sum_{t=1}^\infty t \cdot \mbox{Prob($T=t$)}.
 \end{align} 
 Note now that
 \begin{align*}
\mbox{Prob($T=1$)}&=\sum_{n=1}^N {\bar w}_n^2, \qquad \mbox{Prob($T=2$)}=\left(1-\sum_{n=1}^N {\bar w}_n^2\right)\sum_{n=1}^N {\bar w}_n^2, \quad \mbox{ and }\\
% \quad ... \quad
\mbox{Prob($T=t$)}&=\left(1-\sum_{n=1}^N {\bar w}_n^2\right)^{t-1}\sum_{n=1}^N {\bar w}_n^2.
  \end{align*} 
Thus, replacing into Eq. \eqref{ETformula}, we have 
 \begin{align}\label{ETformula2}
 E[T]&=\sum_{t=1}^\infty t \cdot\left[ \left(1-\sum_{n=1}^N {\bar w}_n^2\right)^{t-1}\sum_{n=1}^N {\bar w}_n^2\right], \\
 &=\left(\sum_{n=1}^N {\bar w}_n^2\right)\left[\sum_{t=1}^\infty t \left(1-\sum_{n=1}^N {\bar w}_n^2\right)^{t-1}\right], \\
 &=\left(\frac{\sum_{n=1}^N {\bar w}_n^2}{1-\sum_{n=1}^N {\bar w}_n^2}\right)\left[\sum_{t=1}^\infty t \left(1-\sum_{n=1}^N {\bar w}_n^2\right)^{t}\right].
 \end{align}
 To simplify the expression above, we can  set $r=1-\sum_{n=1}^N {\bar w}_n^2$, so that we can rewrite it as
  \begin{align}\label{ETformula3}
 E[T]&=\left(\frac{1-r}{r}\right)\left[\sum_{t=1}^\infty t \ r^{t}\right], \\
 &=\left(\frac{1-r}{r}\right) \frac{r}{(1-r)^2}=\frac{1}{1-r}=\frac{1}{\sum_{n=1}^N {\bar w}_n^2}, 
 \end{align}
 where we have used the following equality,
 $$
 \sum_{t=1}^\infty t \ r^{t}=  \frac{r}{(1-r)^2}, \qquad \mbox{when $r\leq 1$},
 $$
 which is a well-known result of power series. 
 
% we draw random pairs of samples with replacement according to the pmf defined  by ${\bar w}_n$, with $n=1,...,N$, the value  $ \frac{1}{\sum_{n=1}^N {\bar w}_n^2}$ is 

%If we have a single maximum, i.e., the maximum value $\max \bar{w}_n$  is reached only with one sample (only for one index $n$),....

\section{Other form for the ESS formula in Eq. \eqref{FirstDef_P} }
\label{App2}

Let us recall $\sum_{n=1}^N {\bar w}_n=1$  so that the arithmetic mean of the normalized weights is always $\mu =\frac{1}{N} \sum_{n=1}^N {\bar w}_n=\frac{1}{N}$.
 Note that the ESS formula  in  Eq. \eqref{FirstDef_P} can 
\begin{align}
\label{FirstDef_Pv2}
\mbox{ESS}_N({\bf {\bar w}})&= \frac{1}{\sum_{n=1}^N {\bar w}_n^2} =\frac{1}{\frac{1}{N}+N\widehat{\sigma}^2}, 
\end{align} 
where $\widehat{\sigma}^2=\frac{1}{N} \sum_{n=1}^N ({\bar w}_n-\mu)^2$ is the variance of the normalized weights. If $\widehat{\sigma}^2=0$, then $\mbox{ESS}_N({\bf {\bar w}})$ reaches the maximum value $N$. We can write:
\begin{align*}
%\label{FirstDef_Pv3}
\frac{1}{\frac{1}{N}+N\widehat{\sigma}^2}&=\frac{1}{\frac{1}{N}+N\left(\frac{1}{N} \sum_{n=1}^N ({\bar w}_n-\mu)^2 \right)}, \\
&=\frac{1}{\frac{1}{N}+N\left( \frac{1}{N} \sum_{n=1}^N {\bar w}_n^2+\frac{1}{N}  \sum_{n=1}^N \mu^2-2\frac{1}{N} \mu \sum_{n=1}^N {\bar w}_n\right)}, \\
&=\frac{1}{\frac{1}{N}+N\left( \frac{1}{N} \sum_{n=1}^N {\bar w}_n^2+\frac{1}{N} N\mu^2-2\mu^2\right)}, \\
&=\frac{1}{\frac{1}{N}+N\left( \frac{1}{N} \sum_{n=1}^N {\bar w}_n^2-\mu^2\right)}, \\
&=\frac{1}{\frac{1}{N}+\left( \sum_{n=1}^N {\bar w}_n^2- N\mu^2\right)}, \\
&=\frac{1}{\frac{1}{N}+ \sum_{n=1}^N {\bar w}_n^2-\frac{1}{N}}=\frac{1}{\sum_{n=1}^N {\bar w}_n^2}.
\end{align*} 
 The equation above proves the equality \eqref{FirstDef_Pv2}.

\end{document}